\documentclass[]{interact}

\usepackage{xcolor}
\usepackage{epstopdf}
\usepackage[caption=false]{subfig}

\usepackage[longnamesfirst,sort]{natbib}
\bibpunct[, ]{(}{)}{;}{a}{,}{,}
\renewcommand\bibfont{\fontsize{10}{12}\selectfont}

\usepackage{algorithm, algorithmic}
\usepackage{enumerate}
\usepackage{natbib}
\usepackage{graphicx,psfrag,epsf}

\def\bbeta{{\bm\beta}}
\def\hbbeta{{\bm{\hat \beta}}}
\def\bx{{\bm x}}
\def\by{{\bm y}}

\def\bw{{\bm w}}

\def\bff{{\bm f}}

\def\bM{{\bm M}}

\def\hbbeta{{\bm{\hat \beta}}}
\def\bc{{\bf c}}
\def\0{{\bf 0}}

\begin{document}


\title{Efficient computational algorithms for approximate optimal designs}

\author{
\name{Jiangtao Duan\textsuperscript{a}, Wei Gao\textsuperscript{a}\thanks{CONTACT Wei Gao. Email: gaow@nenu.edu.cn} , Yanyuan Ma\textsuperscript{b} and Hon Keung Tony Ng\textsuperscript{c}}
\affil{\textsuperscript{a} Key Laboratory for Applied Statistics of MOE, School of Mathematics and\\
      Statistics, Northeast Normal University, Changchun, Jilin 130024, China\\
      \textsuperscript{b} Department of Statistics, The Pennsylvania State University, University Park\\
      \textsuperscript{c} Department of Statistical Science, Southern Methodist University, Dallas
      }
}

\maketitle

\begin{abstract}
In this paper, we propose two simple yet efficient computational algorithms to obtain approximate optimal designs for multi-dimensional linear regression on a large variety of design spaces. We focus on the two commonly used optimal criteria, $D$- and $A$-optimal criteria. For $D$-optimality, we provide an alternative proof for the monotonic convergence for $D$-optimal criterion and propose an efficient computational algorithm to obtain the approximate $D$-optimal design. We further show that the proposed algorithm converges to the $D$-optimal design, and then prove that the approximate $D$-optimal design converges to the continuous $D$-optimal design under certain conditions. For $A$-optimality, we provide an efficient algorithm to obtain approximate $A$-optimal design and conjecture the monotonicity of the proposed algorithm. Numerical comparisons suggest that the proposed algorithms perform well and they are comparable or superior to some existing algorithms.
\end{abstract}

\begin{keywords}
Approximate experimental design \and $D$-optimal \and $A$-optimal \and Regression model
\end{keywords}

\section{Introduction}
Optimal designs are a class of experimental designs that are optimal with respect to some statistical criteria such as minimizing the variance of  best linear unbiased estimators in regression problems and maximizing the amount of information obtained from the experiment.
It is desirable to design experiments that provide more information and reduce the uncertainty relating to the goal of the study.
In regression problems, we model the responses of a random experiment, denoted as $Y_1, \ldots, Y_N$, whose inputs are represented by a vector $\bx_i \in \mathcal{X}$ with respect to some known regression functions $\bff(\bx_1),\ldots,\bff(\bx_N)\in \mathbb{R}^{p}$, i.e.,
$$Y(\bx_i)=\bff^{T}(\bx_i)\bbeta+\varepsilon,i=1,2, \ldots,N,$$
where $\bff(\bx)\in \mathbb{R}^{p}$ is the covariates which is independent variable (regressor) associated with $\bx$, the vector $\bbeta \in \mathbb{R}^{p}$ is $p$-dimensional parameter vector, and $\varepsilon$ is the error term with $E[\varepsilon] = 0$ and  $Var[\varepsilon]=\sigma^{2} > 0$. For different independent trials, the errors are assumed to be uncorrelated and independent, and the inputs $\bx_{i}$ (a candidate set of design points) are chosen by the experimenter in the design space $\mathcal{X}$. We assume that the model is non-singular in the sense that $\{\bff(\bx):\bx\in \mathcal{X}\}$ spans $\mathbb{R}^{p}$. We wish to pick a small subset of the input vector $\bx_i$ such that querying the corresponding responses will lead to a good estimator of the model. In this paper, we discuss the computation of approximate optimal designs for regression models with uncorrelated errors
\citep[see, for example,][]{Fedorov,Harman2018,Atkinson}.

Assume that the inputs $\tilde{\bx}_{i}$, for $i = 1, 2, \ldots, N$ are chosen within a set of distinct points $\bx_{1},\ldots,\bx_{m}$ with integer $p\leq m \ll N$ (but in some special case $m$ can be equal to $N$, see Setting 5), and let $n_{k}$ denote the number of times the particular points $\bx_{k}$ occurs among $\tilde{\bx}_{1},\ldots,\tilde{\bx}_{N}$, and $\tilde{N}=n_{1}+\cdots+n_{m}$ indicates the number of all candidate experiments. The exact experimental design can be summarized by defining a design $\xi_{\tilde{N}}$ as
\begin{equation}
{\xi}_{\tilde{N}} = \left(\begin{array}{ccc}\bx_{1} & \cdots & \bx_{m} \\ \frac{n_1}{\tilde{N}} & \cdots & \frac{n_{m}}{\tilde{N}}\end{array}\right).
\label{DO1}
\end{equation}
In the design ${\xi}_{\tilde{N}}$ in Eq. (\ref{DO1}), the first row gives the points in the design space $\mathcal{X}$ where the input parameters have to be taken and the second row indicates the proportion of the experimental units assigned or the frequencies of the experiments repeated at these points. Strictly speaking, an exact experimental design $\xi$ of size $\tilde{N}$ can be characterized by a probability distribution on $\mathcal{X}$ in which the probability of $\xi_{\tilde{N}}$ occurs at $\bx_{j}$ is $n_{j}/\tilde{N}$.

If $\xi$ follows a continuous probability distribution or a discrete probability distribution
\begin{equation}
{\xi} = \left(\begin{array}{ccc}\bx_{1} & \cdots & \bx_{m} \\ w_{1} & \cdots & w_{m}\end{array}\right),
\label{DO2}
\end{equation}
where $\Pr(\xi=\bx_{j})=w_{j}, j=1,2,\ldots,m$ and $\sum_{j=1}^{m} w_{j} = 1$, then $\xi$ is a continuous design. The goal of the design of experiment theory is then to pick $m$ out of the given $N$ experiments so as to make the most accurate estimate of the parameter $\bbeta$.
For the review and details related to the determination of optimal experimental designs, the readers can refer to \cite{Dette} and the references therein.

Let $\Xi$ be the set of all exact designs or approximate designs (i.e., probability measures) on the design space $\mathcal{X}$; for a given design $\xi$, we denote the information matrix of $\xi$ for the experimental design by
$$
\bM(\xi)=\sum\limits_{\bx\in\mathcal{X}}\xi(\bx)\bff(\bx)\bff^{T}(\bx) \quad \text{or}\quad \bM(\xi)=\int\limits_{\bx\in\mathcal{X}} \xi(\bx)\bff(\bx)\bff^{T}(\bx)d\xi(\bx)
$$
Based on this formulation, an approximate $D$-optimal design for quadratic polynomial was provided by \cite{Chen} when the design space is a circle and \cite{Duan} provided two efficient computational algorithms for optimal continuous experimental designs for linear models.
Under these model assumptions, the Fisher information matrix corresponding to $\bbeta$ is proportional to the information matrix.
Therefore, to obtain the most accurate estimate of certain parameters, we aim to choose the $\xi$ such that $\bM(\xi)$ is maximized according to some criterion. In the following, we will focus on 
a probability measure on $\mathcal{X}$ with support given by the points $\bx_{i}$ and weights $w_{i}$ in Eq. (\ref{DO2}).

In order to obtain the optimal design, a general approach is to consider some generally accepted statistical criteria proposed by \cite{Kiefer1974} namely the $\Phi_q$-criteria.
The $D$-optimality and the $A$-optimality are two of the most commonly used optimality criteria due to their natural statistical interpretations.
It has been shown that the computation of some important prediction-based optimality criteria such as the $I$-optimality criterion \citep{Cook1982,Goos2016} that minimizes the average prediction variance can be converted into the computation of the $A$-optimality  \citep[][Section 10.6]{Atkinson}. In particular, $I$-optimal designs on a finite design space can also be computed using the algorithm developed for $A$-optimality. Thus, in this paper, we focus on the $D$- and $A$-optimal designs where the objective functions are in the form of $\Phi_D(M) = \det(\bM)^{-1}$ and $\Phi_A(M) = tr(\bM)^{-1}$, respectively, for any positive definite matrix $\bM$.
The result by \cite{Welch1982} about the NP-hardness of $D$-optimality is only valid for the exact design problem, while in this paper our aim is to develop efficient computational algorithms for searching the solutions $\xi^{*}$ of the optimization problem $\min \log \Phi_D(M(\xi))$ and $\min \log \Phi_A(M(\xi))$ for $D$-optimality and $A$-optimality, respectively.

Optimal design is at the heart of statistical planning and inference using linear models \citep[see, for example,][]{Box1978}. The theory of optimal designs and the development of numerical computational algorithms for obtaining optimal designs have long been studied in the literature under different scenarios. For instance, \cite{Meyer1995} proposed the coordinate exchange algorithm to construct $D$-optimal and linear-optimal experimental designs for exact design. The algorithm uses a variant of the Gauss-Southwell cyclic coordinate-descent algorithm within the $K$-exchange algorithm to achieve substantive reductions in required computing.
\cite{Gao} developed a general class of the multiplicative algorithms for continuous designs, which can be used to obtain optimal allocation for a general regression model subject to the $D$- and $A$-optimal criteria.
For continuous experimental designs, in general, the continuous factors are generated by the vector $\bff(\bx)$ of linearly independent regular functions where the design points $\bx$ filling the design space $\mathcal{X}$. Then, to choose the optimal design points that maximize the information matrix.

There are many analytical methods for obtaining the approximate optimal designs. \cite{Kiefer1959} introduced the equivalence principle and propose in some cases algorithms to solve the optimization problem. Following the early works of \cite{Karlin1966a}, the case of polynomial regression on a compact interval on $\mathbb{R}$ has been widely studied.
The well-known equivalence theorem of \cite{Kiefer1959} led to the development of a practical algorithm called vertex direction methods (VDMs) for the construction of a $D$-optimal design \citep{Fedorov,Wynn}. They also proved the convergence of the sequence to an optimal (in the appropriate sense) design.
\cite{Silvey} proposed a multiplicative algorithm (MUL) for optimal designs on finite design space, of which the analog in the mixture setting with finite, fixed support is an EM algorithm \citep{Dempster1977}. The VDMs and MUL algorithms all are based on the techniques from differentiable optimization. The general idea is to use directional derivatives to find a direction of improvement, and then employ a line search to determine an optimal step length.
\cite{Yu} proposed the cocktail algorithm, which actually is a mixture of multiplicative, vertex-exchange, and VDM algorithms for $D$-optimum design; it includes a nearest-neighbor exchange strategy that helps to apportion weights between adjacent points and has the property that poor support points are quickly removed from the total support points.
\cite{Harman2018} considered an extension and combination of both the VEM algorithm and the $KL$-exchange algorithm that is used to compute exact designs \citep{Atkinson} and developed the randomized exchange method (REX) for the optimal design problem.

Recent progress in this area has been obtained by employing hybrid methods that alternate between steps of the cocktail algorithm, or by using the randomized exchange method. Following the  work of \cite{Gao}, \cite{Duan} proposed an efficient computational algorithm for computing continuous optimal experimental designs for linear models.

In this paper, we aim to propose a computational algorithm to obtain approximate $D$-optimal designs and a computational algorithm to obtain approximate $A$-optimal designs on any compact design spaces.
This paper is organized as follows. The statistical inference based on a regression model along with the form of an information matrix and variance-covariance matrix for the model parameters are presented in Section 2.
After a review of the $D$- and $A$-optimal criteria, the proposed algorithms and the theoretical results related to the convergence and monotonicity of the proposed algorithms are also presented in Section 2.
Section 3 presents some numerical illustrations with several linear regression models on different types of design spaces which are more general in practical applications for $D$-optimality and $A$-optimality designs.
The proofs of the main results are presented in the Appendix.
\section{Algorithms for Approximate Optimal Designs}
In this section, we introduce the method for searching for optimal designs when regression analysis is used. We focus on the numerical computation of approximate $D$- and $A$-optimal designs. For notation simplicity, we denote $\bff_{i} = \bff(\bx_{i})$, $w_{i}=w(\bx_{i})$, $y_{i}=y(\bx_{i})$ in the following.
Consider the linear regression model
\begin{equation}
Y(\bx)=\bff^{T}(\bx)\bbeta+\varepsilon,\;\bx\in \mathcal{X},
\label{MD}
\end{equation}
where $\bff(\bx)$ is the covariates, $\bbeta $ is a $p$-dimensional parameter vector, $\mathcal{X}$ is the design space and $\varepsilon$ is the error term with mean $0$ and variance $\sigma^2$. When the observations $(\bx, \by = y(\bx))$ are obtained based on the model in Eq. (\ref{MD}), the ordinary least squares estimator of $\bbeta$ can be expressed as
$$
{\hbbeta}=\left[\sum\limits_{i=1}^{N} w_{i}\bff_i\bff_i^{T}\right]^{-1}\sum\limits_{i=1}^{N} w_{i}y_i\bff_i,
$$
where $w_i\geq 0$ is the mass on the point $\bx_i,\;i=1, 2, \ldots,N$ and $\sum_{i=1}^{N} w_i=1$.
The variance of $\hbbeta$ can be obtained as
$$
Var(\hbbeta)=\left[\sum\limits_{i=1}^{N} w_{i}\bff_i\bff_i^{T}\right]^{-1}\sigma^2.
$$

Most of the existing computational algorithms for obtaining optimal designs discretize the underlying continuous space by considering a finite design space $\mathcal{X}=\{\bx_{1}, \cdots, \bx_{N}\} \subset \mathbb{R}^{q}$. These existing algorithms rely on either complex algorithms or advanced mathematical programming solvers. Here, we proposed algorithms that are simple yet effective in obtaining the optimal design for $D$-optimality and $A$-optimality without relying on other complex algorithms or advanced mathematical programming solvers.
\subsection{Algorithms for $D$-optimal Designs}
In an experiment, researchers often wish to estimate the model parameters with the highest precision.
One of the commonly used optimality criteria in experimental design is the $D$-optimal design which maximizes the determinant of the Fisher information matrix, which results in minimum volume for the Wald-type joint confidence region for the model parameters if the variance is known \citep{Gilmour2012}.
Specifically, the $D$-optimal design maximizes the log-determinant of the information matrix, i.e., it minimizes the log-determinant of the asymptotic variance-covariance matrix $Var(\hat\bbeta)$. In other words, the $D$-optimality criterion results in minimizing the generalized variance of the parameter estimates. The $D$-optimal criterion can be described as follows.
\par
\par
\noindent
\textbf{$D$-optimal criterion:}
\begin{equation}
\min\limits_{w_{1},\cdots,w_{N}}\left\{-\log|\sum\limits_{i=1}^{N}w_{i}\bff_{i}\bff_{i}^{T}|:\;\mbox{subject to}\;w_{i}\geq 0\;\mbox{and}\;\sum\limits_{i=1}^{N}w_{i}=1,\; i=1,2,\ldots,N\right\},
\label{DO}
\end{equation}
\par

For $D$-optimality, we can obtain the following result.
\par
\noindent
{\bf Theorem 1.} $\bw^{*}$ is the $D$-optimal solution for Eq. (\ref{DO}) if and only if
$$
\sum\limits_{i=1}^{N} w_{i}\bff^{T}(\bx_{i})\left[\sum\limits_{j=1}^{N} w^{*}_{j}\bff(\bx_{j})\bff^{T}(\bx_{j})\right]^{-1}\bff(\bx_{i})\leq p
$$
for $w_{i}\geq0$ and $\sum\limits_{i=1}^{N} w_{i} = 1,\;i=1,\cdots,N$.

Theorem 1 is a special case of a part of the general equivalence theorem, and the detailed proof is provided in the Appendix.
For the $D$-optimal criterion in Eq. ($\ref{DO}$), we propose the following algorithm to obtain the optimal choice of $\bw^{*}$ based on Theorem 1.
\let\AND\undefined
\begin{algorithm}[H]
	\renewcommand{\algorithmicrequire}{\textbf{Input:}}
	\renewcommand{\algorithmicensure}{\textbf{Output:}}
	\caption{Algorithm for $D$-optimal design}
	\label{alg1}
	\begin{algorithmic}[1]
		\REQUIRE Regressor $\bff(\bx)\in R^{p}$, design space $\mathcal{X}$, stopping parameter $\gamma$ and tuning parameter $\delta$.\\
		\ENSURE Approximate design points $\bx_{s}^{*}$ and corresponding weight $w_{s}^{*}$, $s=1,\cdots,k$.\\
		\STATE Generate random design points $\bx_{i} \in \mathcal{X}$ and corresponding starting weights $w^{(0)}_{i}$, $i=1,\cdots,N$.\\
        \STATE compute the regressors $\bff_{1},\cdots,\bff_{N}\in R^{p}$.\\
       \REPEAT
       \STATE  $$ w^{(h)}_{i}=\frac{w^{(h-1)}_{i}\bff^{T}_{i}D^{(h-1)}\bff_{i}}{p},i=1,\cdots,N,$$
                 where
                 $$ D^{(h-1)}=\left[\sum_{i=1}^{N} w^{(h-1)}_{i}\bff_{i}\bff^{T}_{i}\right]^{-1}.$$\\
       \UNTIL  $\sum_{i=1}^{N} |w_{i}^{(h)} - w_{i}^{(h-1)} | < \gamma$.   \\
       \STATE  find $w_{i}>\delta$ and corresponding design points $\bx_{i}$, $i=1,\cdots,N$, and write them as $\bx^{*}_{s}$, $s=1,\cdots,k$. \\
       \STATE  Let $\bx^{*}_{s}$, $s=1,\cdots,k$ as the new design points and repeat Step 2--5.\\
       \STATE  Output the optimal design points $\bx^{*}_{s}$ and corresponding weights $w^{*}_{s}$, $s=1,\cdots,k$.
	\end{algorithmic}
\end{algorithm}


In Algorithm 1, $\bff(\bx)$ is the functional form of the regressors, the stopping parameter $\gamma$ determines the stopping criteria of the algorithm and the tuning parameter $\delta$ is used to choose the design points with weights which are larger than $\delta$. Note that the choice of $\delta$ must such that $k\geq p$, which is crucial for guaranteeing the non-singularity of the information matrix. In general, the experimenter only needs to set the values of $\gamma$ and $\delta$ to be very small. In Section 3, we set $\delta=0.0001$ and $\gamma=0.0005$.

\noindent
{\bf Remark 1.} Algorithm 1 is efficient even when the sample size $N$ is very large because the iteration process relies on the weight of each sample, but it does not rely on the weights of other samples in a particular iteration. Hence, one can use a parallel strategy to speed up the computations required for Algorithm 1.
Most of the existing algorithms do not share this advantage because the design points interact with each other. For example, the VDM algorithm is based on the differentiable optimization techniques in which the basic idea is to move the current design point $\bx$ to the direction of some other design points while decreasing all components of $\bx$. In addition, the algorithms derived from the VDM algorithm are all depending on the sample size $N$. Therefore, these algorithms will suffer from losing efficiency when the sample size $N$ is large.

In the following, we present proof of the convergence of Algorithm 1 for $D$-optimality. To prove the convergence of the proposed algorithm, we first show that the log-determinant of the information matrix is monotonic. To prove the convergence of the proposed algorithm, we need to add the bounded assumption and require the following two lemmas.\\

\noindent
{\bf Lemma 1.} Let $A(\bx)$ be a nonnegative definite matrix function on $\mathcal{X}$, \linebreak $\bw(\bx)=(w(\bx_1),w(\bx_2),\cdots,w(\bx_N))$ and $\tilde{\bw}(\bx)=(\tilde{w}(\bx_1),\tilde{w}(\bx_2),\cdots,\tilde{w}(\bx_N))$ are two probability vectors in $R^{N}$, and
$$\sum\limits_{i=1}^{N} w(\bx_{i}) A(\bx_{i}) \;\;\mbox{and}\;\;\sum\limits_{i=1}^{N} \tilde{w}(\bx_{i}) A(\bx_{i}) $$
are positive definite matrices. Then,
\begin{eqnarray*}
&      & \log \left|\sum\limits_{i=1}^{N} w(\bx_{i}) A(\bx_{i})\right|-\log \left|\sum\limits_{i=1}^{N} \tilde{w}(\bx) A(\bx_{i}) \right| \\
& \geq & \sum\limits_{i=1}^{N} \tilde{w}(\bx_{i})\mbox{tr}\left\{A(\bx_{i})\left[\sum\limits_{j=1}^{N} \tilde{w}(\bx_{j}) A(\bx_{j}) \right]^{-1}\right\}\log\frac{w(\bx_{i})}{\tilde{w}(\bx_{i})}.
\end{eqnarray*}
\noindent
\textbf{Proof.} Following the proof of Lemma 1 in \cite{Duan} and Lemma 2 in \cite{Gao}, the results in Lemma 1 can be obtained.\\

\noindent
{\bf Lemma 2.} Suppose $\bw(\bx)=(w(\bx_1),w(\bx_2),\cdots,w(\bx_N))$ and \linebreak $\tilde{\bw}(\bx)=(\tilde{w}(\bx_1),\tilde{w}(\bx_2),\cdots,\tilde{w}(\bx_N))$ satisfy $\sum_{i=1}^{N}w(\bx_i)=\sum_{i=1}^{N}\tilde{w}(\bx_i)=1$ are two probability vectors in $R^{N}$, then
$$
\sum\limits_{i=1}^{N}|w(\bx_i)-\tilde{w}(\bx_i)|\leq\left[2\sum\limits_{i=1}^{N} w(\bx_i)\log\frac{w(\bx_i)}{\tilde{w}(\bx_i)}\right]^{1/2}.
$$
\textbf{Proof.} See \cite{Kullback}.  \\

The following theorem shows the convergence of the proposed algorithm for $D$-optimality.\\

\noindent {\bf Theorem 2.} Under the assumption that $\log \left|\sum\limits_{i=1}^{N}\bff(\bx_{i})\bff^{T}(\bx_{i}) \right|$ is bounded, we have
$$
\sum\limits_{i=1}^{N} |w^{(n)}(\bx_{i})-w^{(n-1)}(\bx_{i})|\longrightarrow 0 {\mbox { as }} n \longrightarrow +\infty.
$$
The proof of Theorem 2 is provided in the Appendix.

The algorithm proposed here can be considered as a member of the general class of multiplicative algorithms \citep{Silvey}. Hence, the proposed algorithm shares the simplicity and monotonic convergence property of the class of multiplicative algorithms, and the convergence rate does not depend on $N$ compared to some exact algorithms such as the coordinate-exchange algorithm \citep{Meyer1995}.
Now, we provide a theorem to show that the $D$-optimal approximate design converges to the continuous $D$-optimal design \citep[see, for example,][]{Duan} under certain conditions.
Here, a design is approximate if it is a discrete probability measure and a design is continuous if it is a probability measure with a density with respect to the Lebesgue measure on the observation domain.
Assume that $w^{opt}$ is the $D$-optimal design on the design space $\mathcal{X}$, and it is a continuous probability distribution.\\

\noindent {\bf Theorem 3.} Assume that the random sample points $\bx_i, i=1, 2, \ldots, N$ are generated according to $g(\bx)$ on the design space $\mathcal{X}$ and $g(\bx) > \epsilon > 0$ for all $\bx \in \mathcal{X}$. If $\int_\mathcal{X} \bff^{T}(\bx)\bff(\bx)d\bx$ is bounded on $\mathcal{X}$, then
$$\det\sum\limits_{i=1}^{N} w^*(\bx_i)\bff(\bx_i)\bff^{T}(\bx_i) -\det\int_\mathcal{X} w^{opt}(\bx)\bff(\bx)\bff^{T}(\bx) d\bx\rightarrow 0$$
in probability when $N\to\infty$.\\

\noindent{\bf Proof.} First, we have
\begin{eqnarray*}
\int_\mathcal{X} w^{opt}(\bx) \bff(\bx)\bff^{T}(\bx)d\bx
&=&\int_\mathcal{X} \frac{w^{opt}(\bx)}{g(\bx)} \bff(\bx)\bff^{T}(\bx) g(\bx)d\bx\\
&=&\sum\limits_{i=1}^{N} \frac{w^{opt}(\bx_i)}{Ng(\bx_i)} \bff(\bx_i)\bff^{T}(\bx_i)+\bc_N,
\end{eqnarray*}
where $\bc_N=O_p(N^{-1/2})$. We also have
$0\leq w^{opt}(\bx_i)/\{Ng(\bx_i)\}\leq C/(N\epsilon)<1$
for $N$ sufficiently large.
In addition,
\begin{eqnarray*}
K_N &\equiv& \sum\limits_{i=1}^{N} w^{opt}(\bx_i)/\{Ng(\bx_i)\}\\
&=&\int_\mathcal{X} w^{opt}/g(\bx) g(\bx)d\bx+d_N = 1+d_N,
\end{eqnarray*}
where $d_N=O_p(N^{-1/2})$. Note that $w^{opt}(\bx_i)/\{NK_Ng(\bx_i)\}, i=1, 2, \ldots, N$ are valid weights.

Write $\sum\limits_{i=1}^{N} w^*(\bx_i)\bff(\bx_i)\bff^{T}(\bx_i) = \int_\mathcal{X} w_N(\bx)\bff(\bx)\bff^{T}(\bx) d\mu(\bx)$, where
$w_N(\bx)=w^*(\bx_i)$ if $\bx\in\{\bx_1, \bx_2, \ldots, \bx_N\}$ and $w_N(\bx)=0$
otherwise. Then, because $w^{opt}$ is optimal on $\mathcal{X}$ among all possible bounded probability density functions and probability mass functions on $\mathcal{X}$, we can obtain
\begin{equation}
\det \int_\mathcal{X} w^{opt}(\bx)\bff(\bx)\bff^{T}(\bx) d\bx
\ge \det\sum\limits_{i=1}^{N} w^*(\bx_i)\bff(\bx_i)\bff^{T}(\bx_i).
\label{eq:relation1}
\end{equation}
Given $\bx_i, i = 1, 2, \ldots, N$, $w^*(\bx_i)$ is
optimal, we have
\begin{equation}
\det \sum\limits_{i=1}^{N} \frac{w^{opt}(\bx_i)}{NK_Ng(\bx_i)} \bff(\bx_i)\bff^{T}(\bx_i)
\le \det \sum\limits_{i=1}^{N} w^*(\bx_i)\bff(\bx_i)\bff^{T}(\bx_i).
\label{eq:relation2}
\end{equation}
Combining Eqs. (\ref{eq:relation1}) and (\ref{eq:relation2}),
we have
\begin{eqnarray}
\nonumber
\det \int_\mathcal{X} w^{opt}(\bx)\bff(\bx)\bff^{T}(\bx) d\bx &\geq& \det \sum\limits_{i=1}^{N} w^*(\bx_i)\bff(\bx_i)\bff^{T}(\bx_i)\\
\nonumber
&\geq& \det \frac{1}{K_N}\sum\limits_{i=1}^{N} \frac{w^{opt}(\bx_i)}{Ng(\bx_i)} \bff(\bx_i)\bff^{T}(\bx_i)\\
&=&\det \frac{1}{1+d_N}\left\{ \int_\mathcal{X} w^{opt}(\bx)\bff(\bx)\bff^{T}(\bx) d\bx-\bc_N\right\}.
\label{eq:combined}
\end{eqnarray}
Because $\bc_N\rightarrow\0$ and $d_N\rightarrow0$ in probability when $N\rightarrow\infty$, Eq. (\ref{eq:combined}) leads to
\begin{eqnarray*}
&&\det\int_\mathcal{X} w_N(\bx)\bff(\bx)\bff^{T}(\bx) d\mu(\bx)
-\det\int_\mathcal{X} w^{opt}(\bx)\bff(\bx)\bff^{T}(\bx) d\bx\\
&=&\det\sum\limits_{i=1}^{N} w^*(\bx_i)\bff(\bx_i)\bff^{T}(\bx_i)
-\det\int_\mathcal{X} w^{opt}(\bx)\bff(\bx)\bff^{T}(\bx) d\bx\rightarrow0
\end{eqnarray*}
in probability when $N\rightarrow\infty$. \\

Theorem 3 verifies that the approximate optimal design obtained from the proposed algorithm will eventually converge to the continuous optimal design at the speed of $\sqrt{N}$. Theorem 3 also guarantees the convergence of the proposed algorithm.
\subsection{Algorithms for A-optimal Designs}
Another commonly used optimality criterion is the $A$-optimality criterion that minimizes the trace of the variance-covariance matrix of the maximum likelihood estimates (MLEs). The $A$-optimality criterion provides an overall measure of the variations in the model parameter estimates.
The objective function being minimized in the $A$-optimal design is described as follows.  \\
\textbf{$A$-optimal criterion:}
\begin{equation}
\min\limits_{w_1,\cdots,w_{N}} \left\{ {\mbox {tr}} \left( \left[\sum\limits_{i=1}^{N} w_{i}\bff_{i}\bff^{T}_{i} \right]^{-1} \right), \mbox{subject to}\;w_{i}\geq 0\;\mbox{and}\;\sum\limits_{i=1}^{N} w_{i}=1,\; i=1,\cdots,N \right\}.
\label{AOP}
\end{equation}

For $A$-optimality, we can obtain the following theorem. For simplicity, we denote $ I_\mathcal{X}(\bw, \bff)  = \sum\limits_{i=1}^{N} w(\bx_{i})\bff_{i}(\bx_{i})\bff_{i}^{T}(\bx_{i})$ and $I_\mathcal{X}(\bw^{*}, \bff)=\sum\limits_{i=1}^{N} w^{*}(\bx_{i})\bff_{i}(\bx_{i})\bff^{T}(\bx_{i})$.

\noindent
{\bf Theorem 4.} $\bw^{*}(\bx)=(w^{*}(\bx_1),\cdots,w^{*}(\bx_N))$ is the $A$-optimal solution for Eq. (\ref{AOP}) if and only if
$$
\frac{ {\mbox {tr}}(\left[I_\mathcal{X}(\bw^{*}, \bff) \right]^{-1} I_\mathcal{X}(\bw, \bff) \left[ I_\mathcal{X}(\bw^{*}, \bff) \right]^{-1})} {{\mbox {tr}}(\left[I_\mathcal{X}(\bw^{*}, \bff)\right]^{-1})} \leq 1,
$$
for $ w(\bx_i)\geq0, \; i=1,2,\ldots,N$ and $\sum\limits_{i=1}^{N} w(\bx_i)=1$.

The proof of Theorem 4 is provided in the Appendix. A similar approach to Theorem 3 can be used to show that the approximate $A$-optimal design converges to the continuous $A$-optimal design. Based on Theorem 4, for the $A$-optimal criteria in Eq. (\ref{AOP}), the following algorithm is proposed to obtain the $A$-optimal design $\bw^{*}(\bx)$.\\

\let\AND\undefined
\begin{algorithm}[H]
	\renewcommand{\algorithmicrequire}{\textbf{Input:}}
	\renewcommand{\algorithmicensure}{\textbf{Output:}}
	\caption{Algorithm for $A$-optimal design}
	\label{alg2}
	\begin{algorithmic}[1]
		\REQUIRE Regressor $\bff(\bx)\in R^{m}$, design space $\mathcal{X}$, stopping parameter $\gamma$ and tuning parameter $\delta$.\\
		\ENSURE Approximate design points $\bx_{s}^{*}$ and corresponding weight $w_{s}^{*}$, $s=1,2,\ldots,k$.\\
		\STATE Generate a random design points $\bx_{i} \in \mathcal{X}$ and corresponding starting weights $w^{(0)}_{i}$, $i=1,2,\ldots,N$.\\
        \STATE compute the regressors $\bff_{1},\bff_{2},\ldots,\bff_{N}\in R^{p}$.\\
       \REPEAT
       \STATE    {$$ w^{(h)}_{i}= w^{(h-1)}_{i} \left[\frac{(p-1)}{p} \frac{tr(D^{(h-1)}\bff_{i}\bff_{i}^{T} D^{(h-1)} )}{tr(D^{(h-1)})} +\frac{1}{p} \right], i=1,2,\ldots,N,$$}
                 where
                 $$ D^{(h-1)}=\left[\sum_{i=1}^{N} w^{(h-1)}_{i}\bff_{i}\bff^{T}_{i}\right]^{-1}.$$\\
       \UNTIL  $\sum_{i=1}^{N} |w_{i}^{(h)} - w_{i}^{(h-1)} | < \gamma$.   \\
       \STATE  find $w_{i}>\delta$ and corresponding design points $\bx_{i}$, $i=1,2,\ldots,N$, and write them as $\bx^{*}_{s}$, $s=1,2,\ldots,k$. \\
       \STATE  Let $\bx^{*}_{s}$, $s=1,2,\ldots,k$ as the new design points and repeat Step 2--5.\\
       \STATE  Output the optimal design points $\bx^{*}_{s}$ and corresponding weights $w^{*}_{s}$, $s=1,\cdots,k$.
	\end{algorithmic}
\end{algorithm}

In Algorithm 2, $\bff(\bx)$ is the functional form of the regressors, the stopping parameter $\gamma$ determines the stopping criteria of the algorithm and the tuning parameter $\delta$ is to choose the optimal design points with the weight larger than $\delta$. For the same reason as described in Remark 1, a parallel strategy can be used to speed up the computations required for Algorithm 2.

For $A$-optimality, the proposed algorithm provides a convergence solution that is robust to the initial value because the algorithm does not depend on the initial value.
We have attempted to develop the theoretical justification of the convergence of the proposed computational algorithm for $A$-optimality, however, the mathematical justification is not available. Instead, we provide some simulation and numerical results to support the validity and reliability of the proposed algorithm in the subsequent section.
%
Here, we conjecture the monotonic convergence of the algorithm for $A$-optimality based on the extensive simulation and numerical results.
\section{Numerical Illustrations}
In this section, we consider five different settings \citep{Castro2017} to evaluate the performance of the proposed algorithms for obtaining approximate $D$- and $A$-optimal designs.
We use these numerical examples to illustrate that the proposed algorithms can efficiently identify the optimal design.

For comparative purposes, we also apply the randomized exchange (REX) algorithm  proposed by \cite{Harman2018}, the cocktail (CO) algorithm proposed by \cite{Yu}, the vertex direction (VDM) algorithm proposed by \cite{Fedorov} and \cite{Wynn}, and the multiplicative (MUL) algorithm proposed by \cite{Silvey} for computing the $D$- and $A$-optimal designs.
Since the cocktail algorithm is for $D$-optimal only \citep{Yu}, hence, the cocktail algorithm is not applied to obtain the $A$-optimal design.
\subsection{Setting 1: Two-dimensional design space in a square with $p = 6$}
In Setting 1, we consider the model
$$
y(\bx)=\bbeta^{T}\bff(\bx),
$$
where $\bbeta=(\beta_0,\beta_1,\beta_2,\beta_3,\beta_4,\beta_5)^{T}$, $\bff(\bx)=(1,x_1,x_2,x_1^2,x_1x_2,x_2^2)^{T}$ and $\bx\in \mathcal{X}=[-1,1] \times [-1,1]$.
The optimal design points and their corresponding weights obtained from Algorithm \ref{alg1}, Algorithm \ref{alg2} and the REX, CO, VEM and MUL algorithms for $D$- and $A$-optimal designs are presented in Tables \ref{set1_DOP} and \ref{set1_AOP}, respectively. The optimal values of the corresponding objective functions are also presented in Tables \ref{set1_DOP} and \ref{set1_AOP}.

From Tables \ref{set1_DOP} and \ref{set1_AOP}, we observe that the proposed algorithms identify the same optimal design points as the REX, CO, VDM, and MUL algorithms. The weights for optimal points obtained from the REX, CO algorithms and the proposed algorithm are very close, and there is no significant difference between these three algorithms according to the values of the $D$-optimality objective function presented in Table \ref{set1_DOP}. Furthermore, the performance of the VDM and MUL algorithms are not as good as the REX, CO and proposed algorithms. Similar results and conclusions can be observed from Table \ref{set1_AOP} for $A$-optimal designs.
\subsection{Setting 2: Two-dimensional design space in a circle with $p = 6$}
In Setting 2, we consider the model
\begin{eqnarray*}
y(\bx)=\bbeta^{T}\bff(\bx),
\end{eqnarray*}
where $\bbeta=(\beta_0,\beta_1,\beta_2,\beta_3,\beta_4,\beta_5)^{T}$, $\bff(x)=(1,x_1,x_2,x_1^2,x_1x_2,x_2^2)^{T}$ and $\bw\in \mathcal{X}=\{(x_1,x_2)^{T}: \;x_1^2+x_2^2\leq 1\}$.

\begin{figure}[H]
  \centering
  \includegraphics[width=0.9\textwidth]{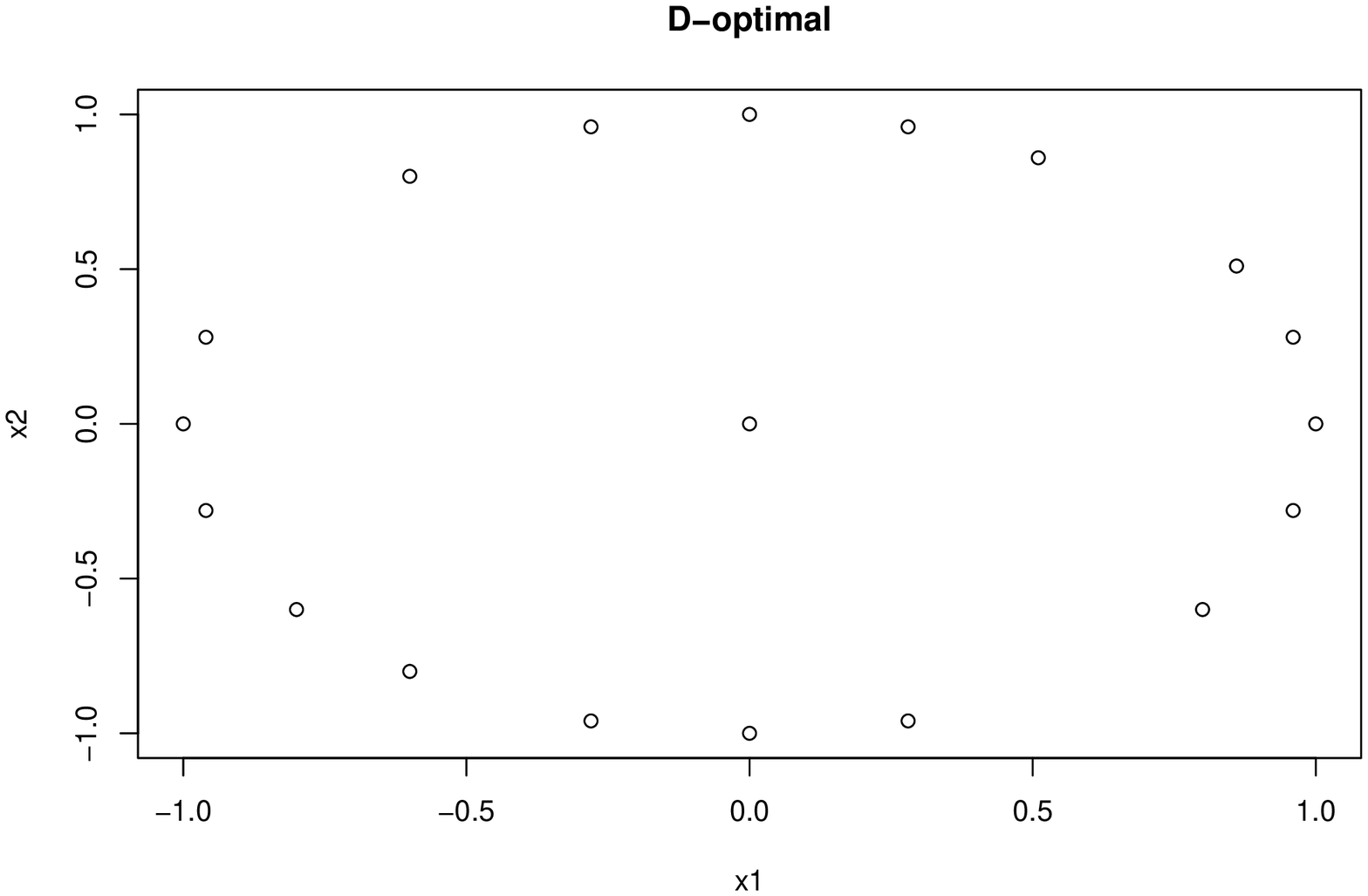}\\
  \caption{$D$-optimal design for Setting 2}\label{set2}
\end{figure}

Figure \ref{set2} presents the $D$-optimal design points in which the weight for the center of the unit circle (0, 0) is 1/6, and the other optimal design points are uniformly distributed on the ring (the vertices of a regular $s$-sided polygon in the circle) with a combined weight $5/6$ in theory. For a more detailed analysis of this setting, one can refer to \cite{Duan}. The optimal design points and their corresponding weights obtained from Algorithm \ref{alg1}, Algorithm \ref{alg2} and the REX, CO, VEM and MUL algorithms for $D$- and $A$-optimal designs are presented in Tables \ref{set2_DOP} and \ref{set2_AOP}. The optimal values of the corresponding objective functions are also presented in Tables \ref{set2_DOP} and \ref{set2_AOP}, respectively.
Note that different algorithms may have different design points because the optimal design points are distributed uniformly on a circle with almost equal weight. Hence, the points may locate at different locations that are symmetrical about the center of the circle.

From Table \ref{set2_DOP}, the REX and CO algorithms, and the proposed algorithm have the same weight as the theoretical value for the center point $(0, 0)$ and the other optimal design points are evenly distributed on the ring with total weight $5/6$, while the VDM and MUL algorithms have distributed the weights for all design points including the points in the center and on the ring.
By comparing the values of the objective function $-\log(|\Sigma(\bx)|)$ of different algorithms, we can see that the proposed algorithm and the REX algorithm give the same value of the objective function, which is a better value compared to the values obtained from other algorithms. Once again, similar results and conclusions can be observed from Table \ref{set2_AOP} for $A$-optimal designs.

\subsection{Setting 3: Two-dimensional Wynn's polygon design space with $p = 6$}
In Setting 3, we consider a two-dimensional irregular design space called Wynn's polygon as
$$
y(\bx)=\bbeta^{T}\bff(\bx),
$$
where $\bbeta=(\beta_0,\beta_1,\beta_2,\beta_3,\beta_4,\beta_5)^{T}$, $\bff(\bx)=(1,x_1,x_2,x_1^2,x_1x_2,x_2^2)^{T}$ and $\bx\in \mathcal{X}=\{(x_1,x_2)^{T}: \;x_1,x_2\geq -\frac{1}{4}\sqrt{2},x_1\leq(x_2+\sqrt{2}),x_2\leq \frac{1}{3}(x_1+\sqrt{2}),x_1^2+x_2^2\leq 1)\}$.\\

\begin{figure}[H]
  \centering
  \includegraphics[width=0.9\textwidth]{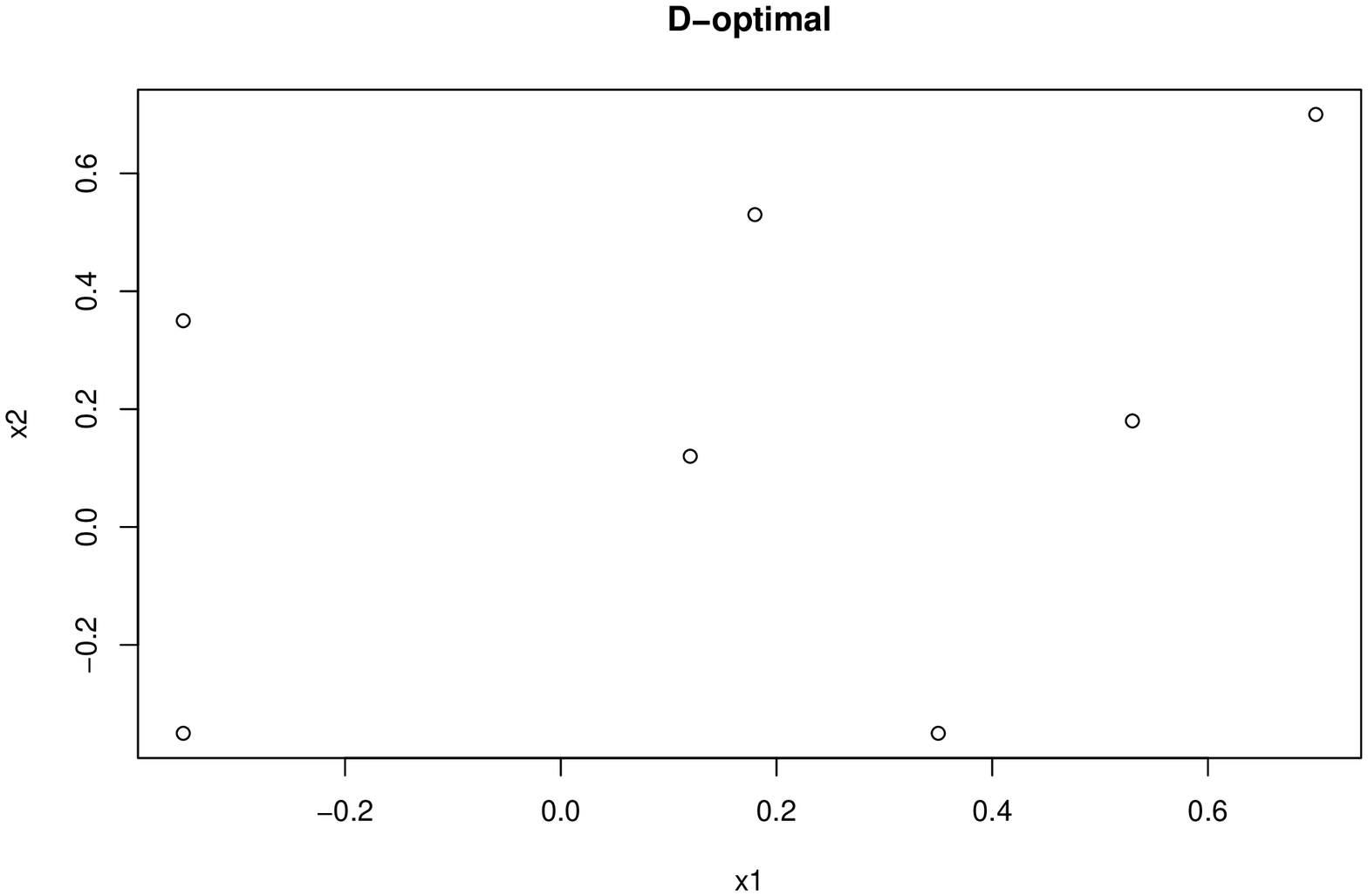}\\
  \caption{$D$-optimal design for Setting 3}\label{set3}
\end{figure}

Figure \ref{set3} shows the $D$-optimal design points and Table \ref{set3_DOP} shows the design points and corresponding weights by different algorithms. From Table \ref{set3_DOP}, we can see that the proposed algorithm, the REX and CO algorithm can locate the optimal design points with the theoretical weights while the  VDM and multiplicative algorithm fail to do so.
For $A$-optimal, the design points and corresponding weights by different algorithms are presented in Table \ref{set3_AOP}. From Table \ref{set3_AOP}, we observe that the proposed algorithm and the REX algorithm provide the same results, and they are superior to other methods.

\subsection{Setting 4: Three-dimensional design space in a cube with $p = 10$}
In Setting 4, we consider the model
$$y(\bx)=\bbeta^{T}\bff(\bx),$$
where $\bbeta=(\beta_0,\beta_1,\beta_2,\ldots, \beta_8,\beta_9)^{T}$, $\bff(\bx)= (1,x_1,x_2,x_3,x_1^2,x_1x_2,x_1x_3,x_2^2,x_2x_3,x_3^2)^{T}$ and $\bx\in \mathcal{X}=[-1,1]\times[-1,1]\times[-1,1]$.

Based on the algorithms considered here, we obtain the same 27 $D$- and $A$-optimal design points with different weights. The weights obtained from different algorithms with the corresponding values of the objective function are presented in Tables \ref{set4_DOP} and \ref{set4_AOP} for the $D$- and $A$-optimality, respectively.
From Tables  \ref{set4_DOP} and \ref{set4_AOP}, we observe that the results obtained from the proposed algorithms are very close to the theoretical values \cite{Duan}. Note that although the weights of the proposed optimal design are different from the weights of continuous optimal design provided by \cite{Atkinson}, our optimal design has a smaller $D$-optimal value. Therefore, according to the definition of $D$-optimal criterion, the proposed optimal design should be better.
In this numerical study, we observe that the optimal design points obtained by the existing algorithms are unstable in the sense that the optimal design points may not be unique in multiple runs of the algorithms. This may lead to uncertainty in practical applications. Moreover, we found that the REX, VDM, and MUL algorithms cannot always get all the optimal design points. This may lead to severe problems in some critical experiments. For example, in pharmaceutical or chemical experiments, ignoring some design points may lead to severe consequences.
For illustrative purpose, the $D$-optimal value for the proposed, REX and CO algorithms presented in Tables \ref{set4_DOP} and \ref{set4_AOP} are closest to the theoretical values.

In fact, the REX method may fail when the number of candidate points is small and the number of parameters is large. For example, in Setting 4, if the number of candidate points has 27 design points, the REX method sometimes has a singularity of the $M$ matrix during the calculation process, which causes the failure of obtaining the optimal design. Similarly, in Setting 1, when there are only 7 candidate points, we found that the REX method may also fail sometimes. In contrast, the proposed method is feasible and stable in obtaining the optimal design points even when the number of candidate points is small.
This is a significant advantage of the proposed method in practical application because there are many scenarios that the number of candidate points is small due to high cost or environmental factors.



\subsection{Setting 5: Three-dimensional design space in a sphere with $p = 9$}
In Setting 5, we consider the model
$$y(\bx)=\bbeta^{T}\bff(\bx),$$
where $\bbeta=(\beta_1,\beta_2, \ldots, \beta_8,\beta_9)^{T}$, $\bff(\bx)=(x_1,x_2,x_3,x_1^2,x_1x_2,x_1x_3,x_2^2,x_2x_3,x_3^2)^{T}$ and $\bx\in \mathcal{X}=\{(x_1,x_2,x_3)^{T}: \;x_1^2+x_2^2+x_3^2=1\}$.

\begin{figure}[H]
  \centering
  \includegraphics[width=0.8\textwidth]{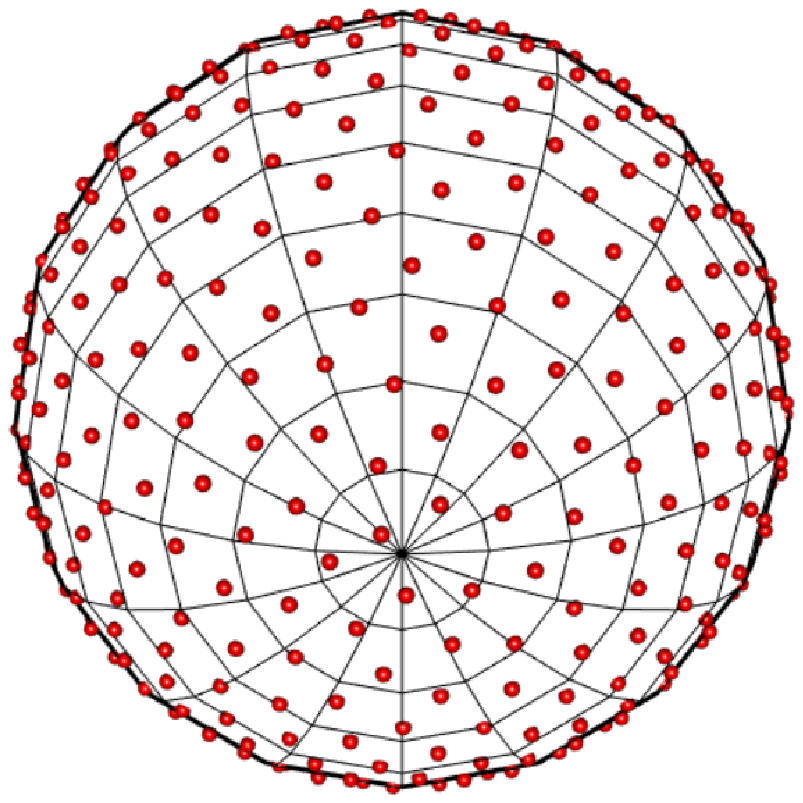}
  \caption{$D$-optimal design for Setting 5}\label{set5_1}
\end{figure}

\begin{figure}[H]
  \centering
  \includegraphics[width=0.8\textwidth]{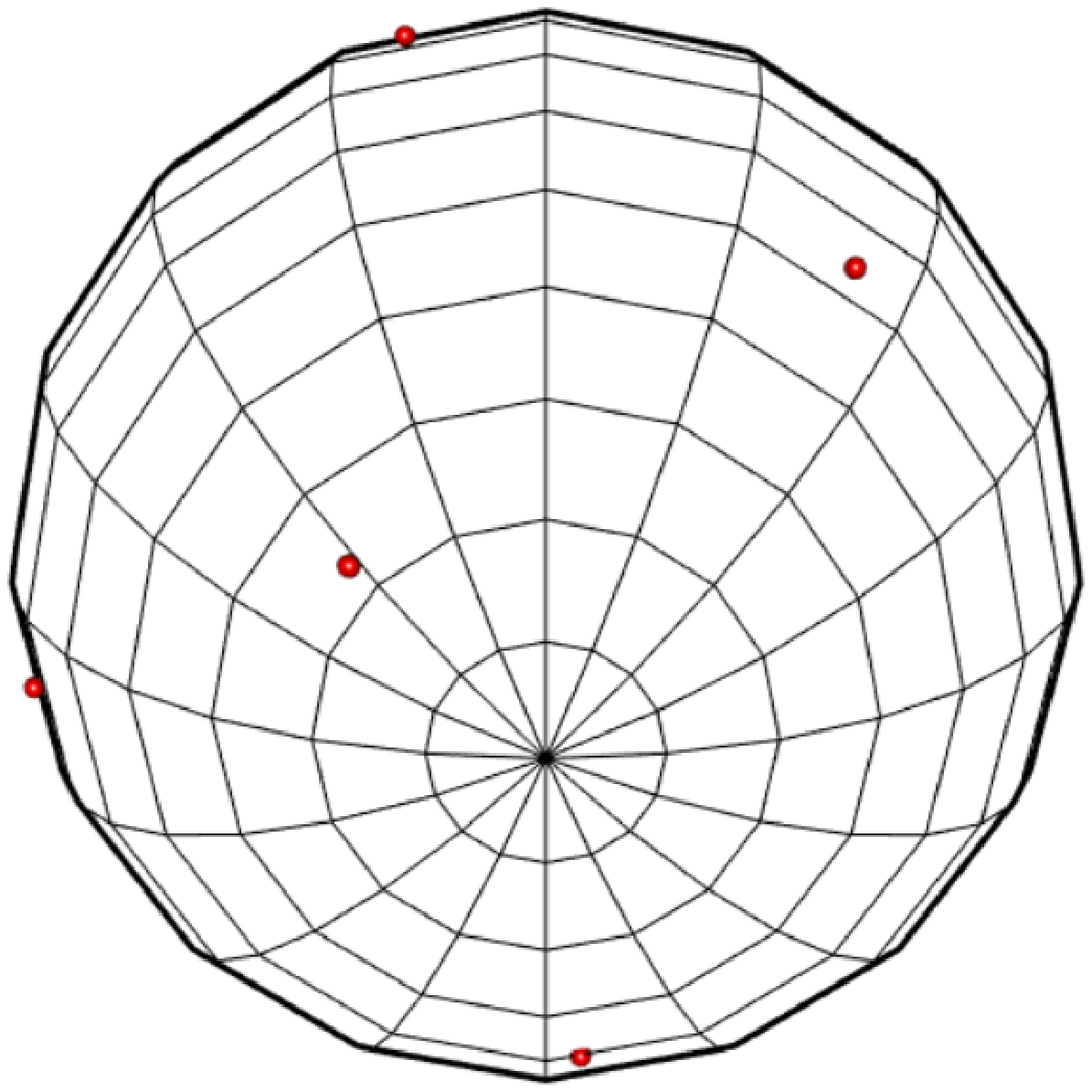}\\
  \caption{$D$-optimal design for Setting 5}\label{set5_2}
\end{figure}

\begin{figure}[H]
  \centering
  \includegraphics[width=0.8\textwidth]{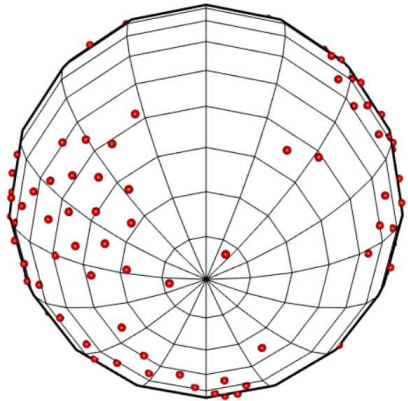}\\
  \caption{$D$-optimal design for Setting 5 by REX method}\label{set5_3}
\end{figure}

For this setting, we use the Fibonacci numbers on the 3-dimensional unit sphere as the initial design points for the algorithms considered here. For more details related to this setting, the reader can refer to \cite{Castro2017}.
Regression problems with a unit sphere design space have many applications in astrophysics, gravity induction, geophysics, climate laws, and global navigation, because there are countless signals on the surface of the earth, and satellite signals also affect our daily lives. Another important application of regression with a unit sphere design space is three-dimensional human faces recognition with sparse spherical representation in authentication and surveillance.
Based on this setting, we find that every design point on the unit sphere can be considered as an optimal design point. By using the proposed algorithms, we obtain all the design points with equal weights.
Figures \ref{set5_1} and \ref{set5_2} display the design points of the $D$-optimal design when the number of supporting points are $500$ and $10$, respectively.
However, when applying the REX, VDM, MUL, and CO algorithms, only a smaller number of design points are identified. For instance, the REX algorithm gives only 128 points as the optimal design points when we use 5000 supporting points in the three-dimensional unit sphere. Moreover, the weights assigned to these 128 design points are not equal based on the REX algorithm.
Figure \ref{set5_3} presents the $D$-optimal design obtained from the REX algorithm when the number of supporting points is 500. We also present the values of the $D$- and $A$-optimality objective functions in Table \ref{set5_DEFF}. From Table \ref{set5_DEFF}, we observe that the $D$-optimal and $A$-optimal values of the proposed algorithm are smaller than the other algorithms considered here. Thus, the proposed method performs well in this case.

\begin{figure}[H]
  \centering
  \includegraphics[width=0.6\textwidth]{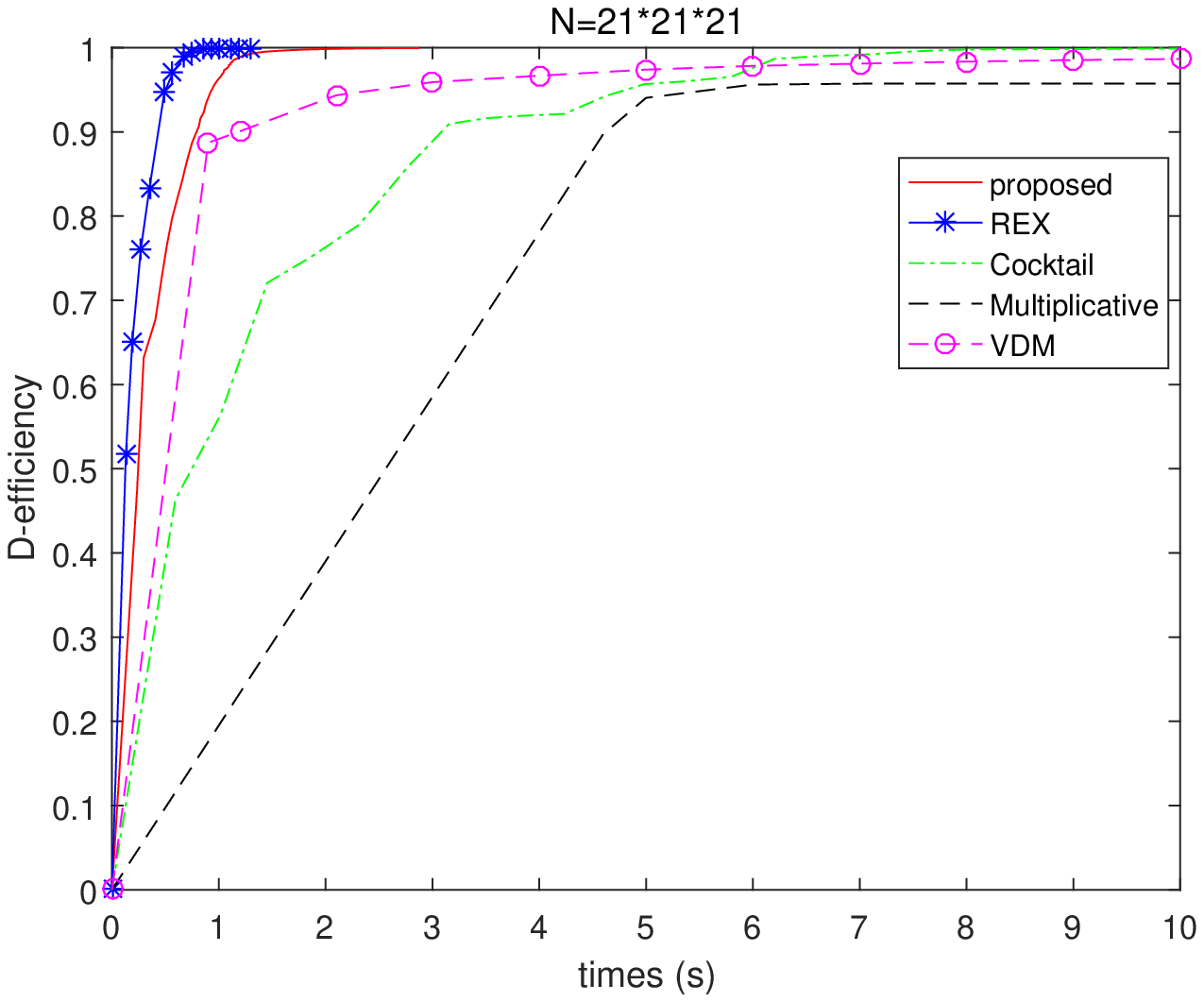}\\
  \caption{The $D$-efficiency for Setting 4 when $N=21\times 21 \times21$. The vertical axis denotes the $D$-efficiency. The horizontal axis corresponds to the computation time (in seconds).}\label{DOP_N21}
\end{figure}

\begin{figure}[H]
  \centering
  \includegraphics[width=0.6\textwidth]{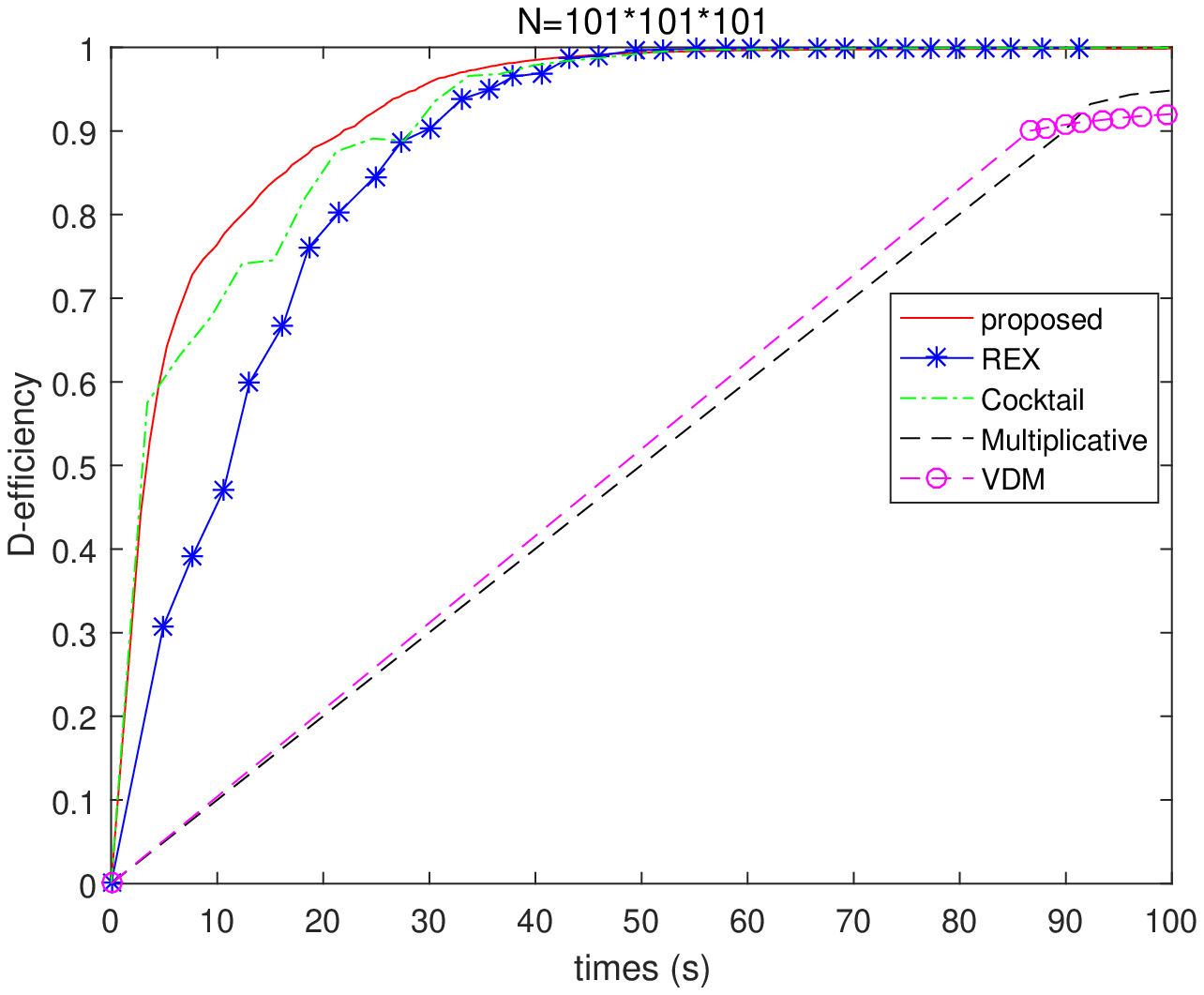}\\
  \caption{The $D$-efficiency for Setting 4 when $N=101\times 101 \times101$. The vertical axis denotes the $D$-efficiency. The horizontal axis corresponds to the computation time (in seconds).}\label{DOP_N101}
\end{figure}

\begin{figure}[H]
  \centering
  \includegraphics[width=0.6\textwidth]{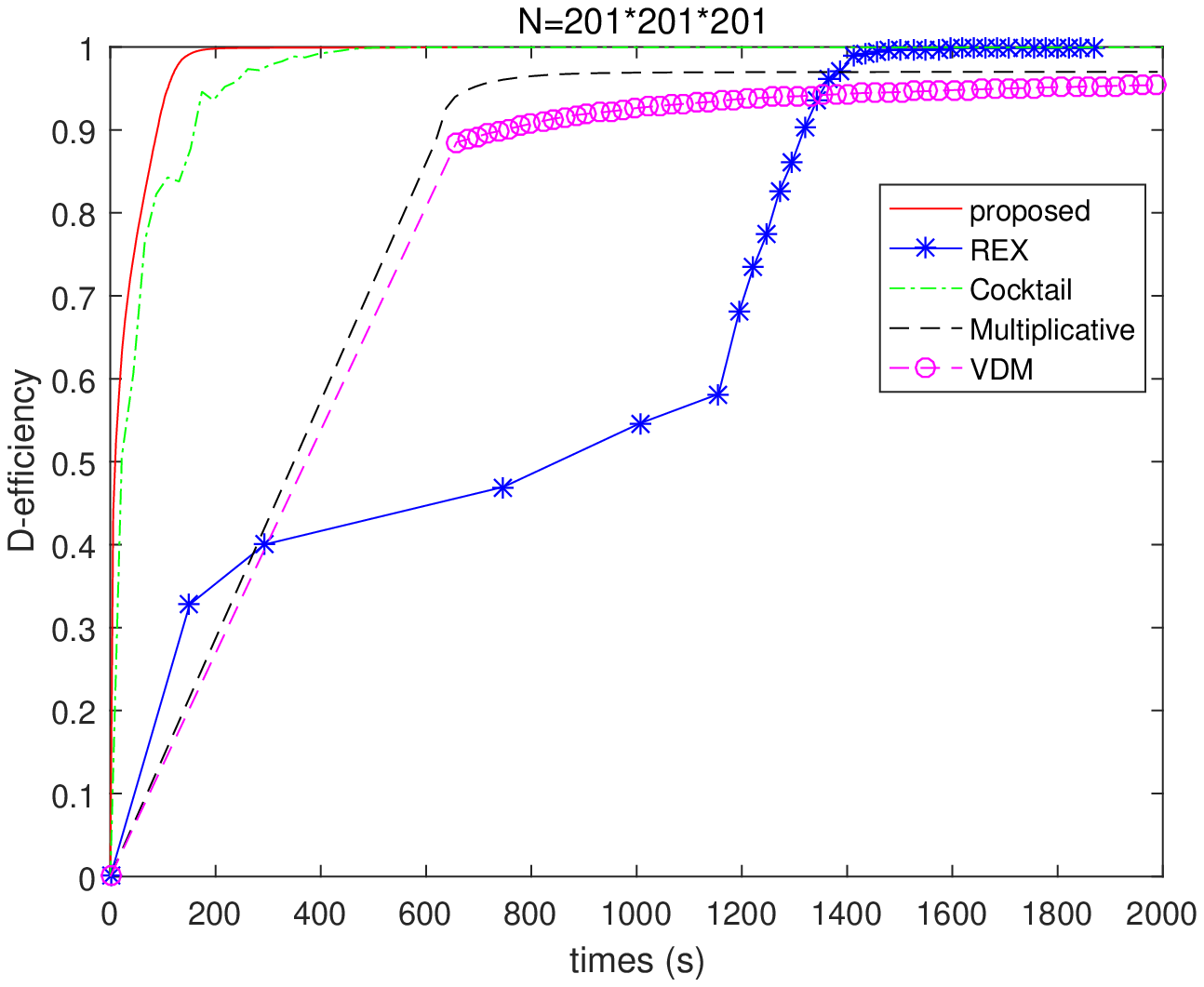}\\
  \caption{The $D$-efficiency for Setting 4 when $N=201\times 201 \times201$. The vertical axis denotes the $D$-efficiency. The horizontal axis corresponds to the computation time (in seconds).}\label{DOP_N201}
\end{figure}

\begin{figure}[H]
  \centering
  \includegraphics[width=0.6\textwidth]{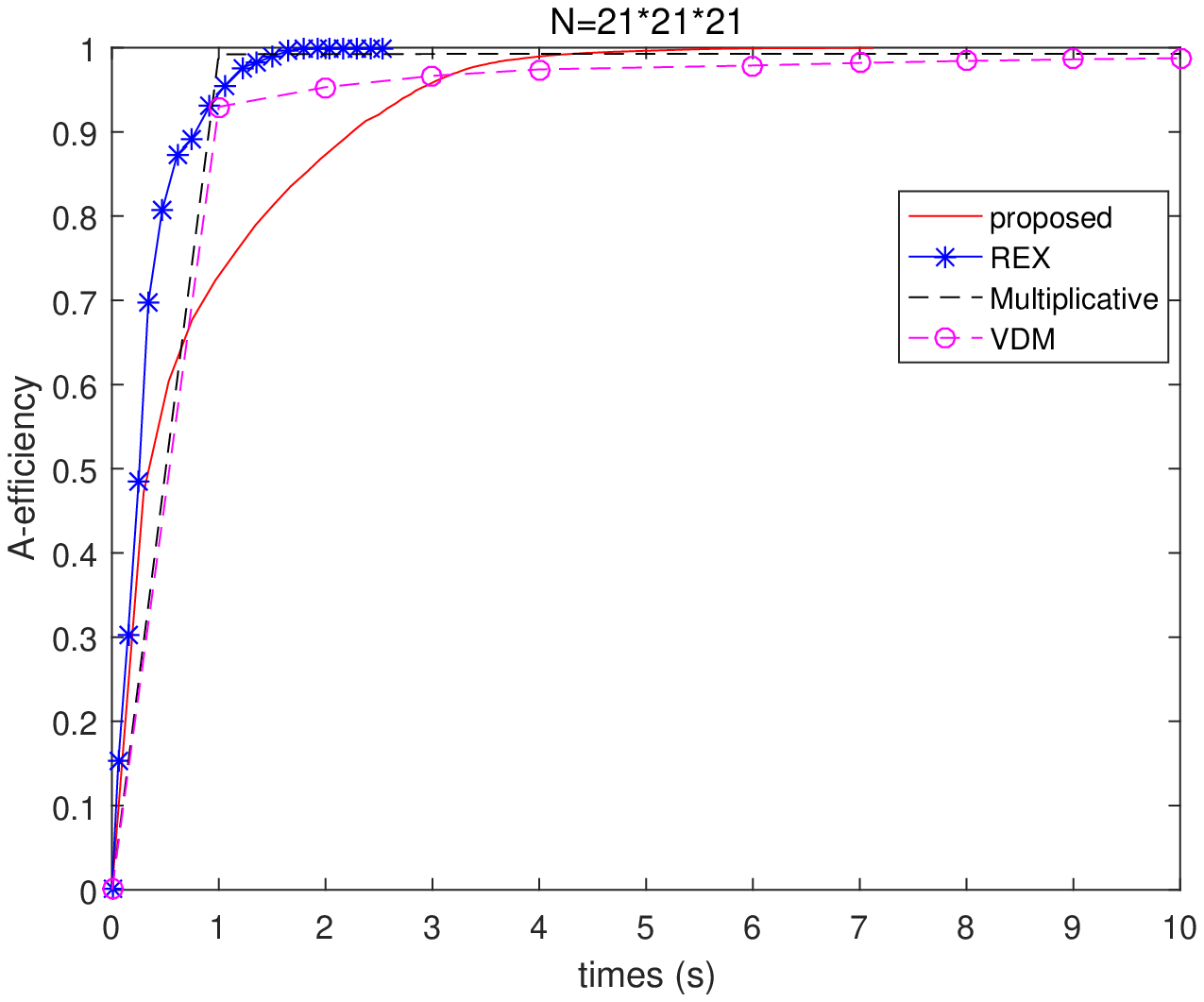}\\
  \caption{The $A$-efficiency for Setting 4 when $N=21\times 21 \times21$. The vertical axis denotes the $A$-efficiency. The horizontal axis corresponds to the computation time (in seconds).}\label{AOP_N21}
\end{figure}

\begin{figure}[H]
  \centering
  \includegraphics[width=0.6\textwidth]{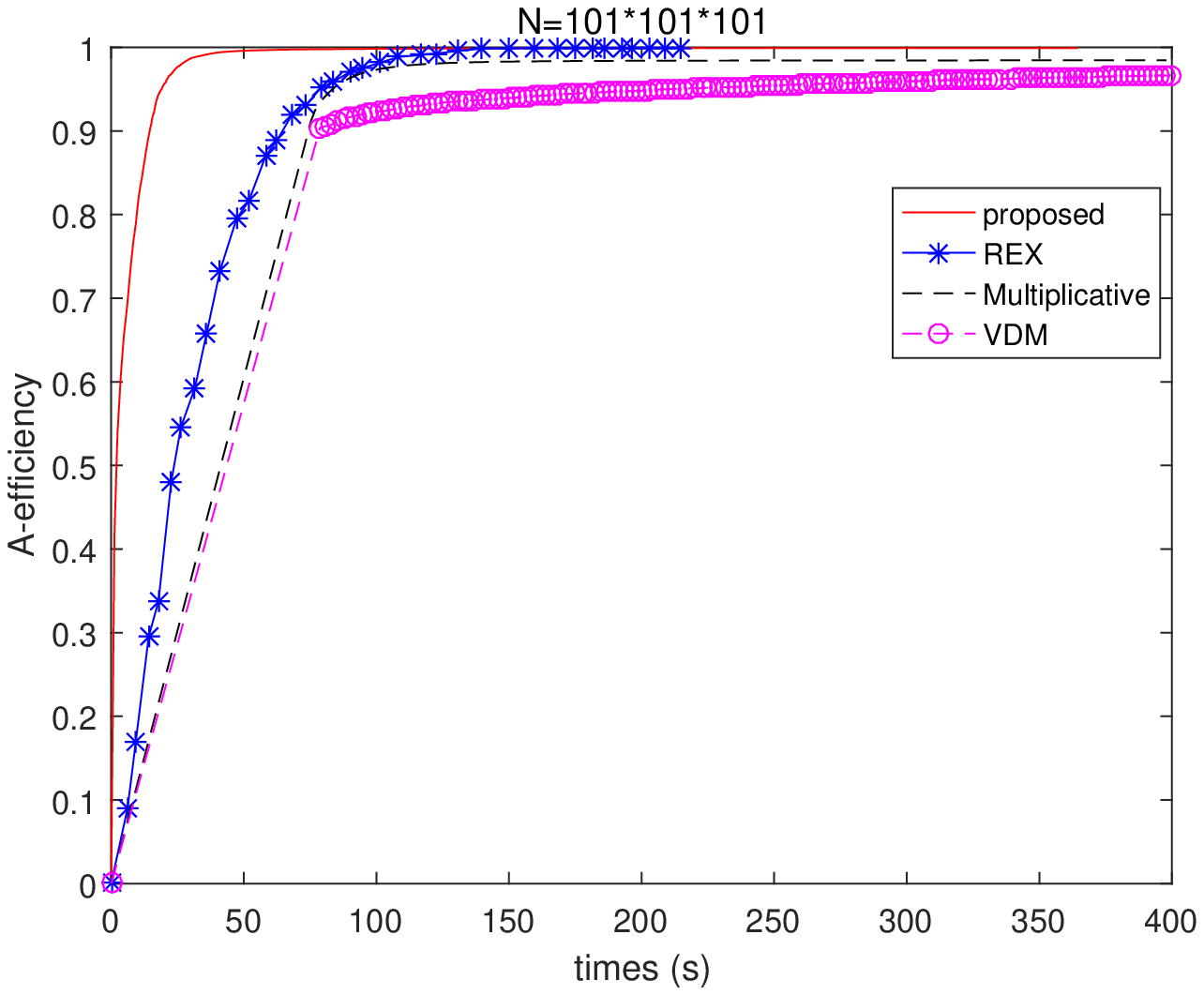}\\
  \caption{The $A$-efficiency for Setting 4 when $N=101\times 101 \times101$. The vertical axis denotes the $A$-efficiency. The horizontal axis corresponds to the computation time (in seconds).}\label{AOP_N101}
\end{figure}

\begin{figure}[H]
  \centering
  \includegraphics[width=0.6\textwidth]{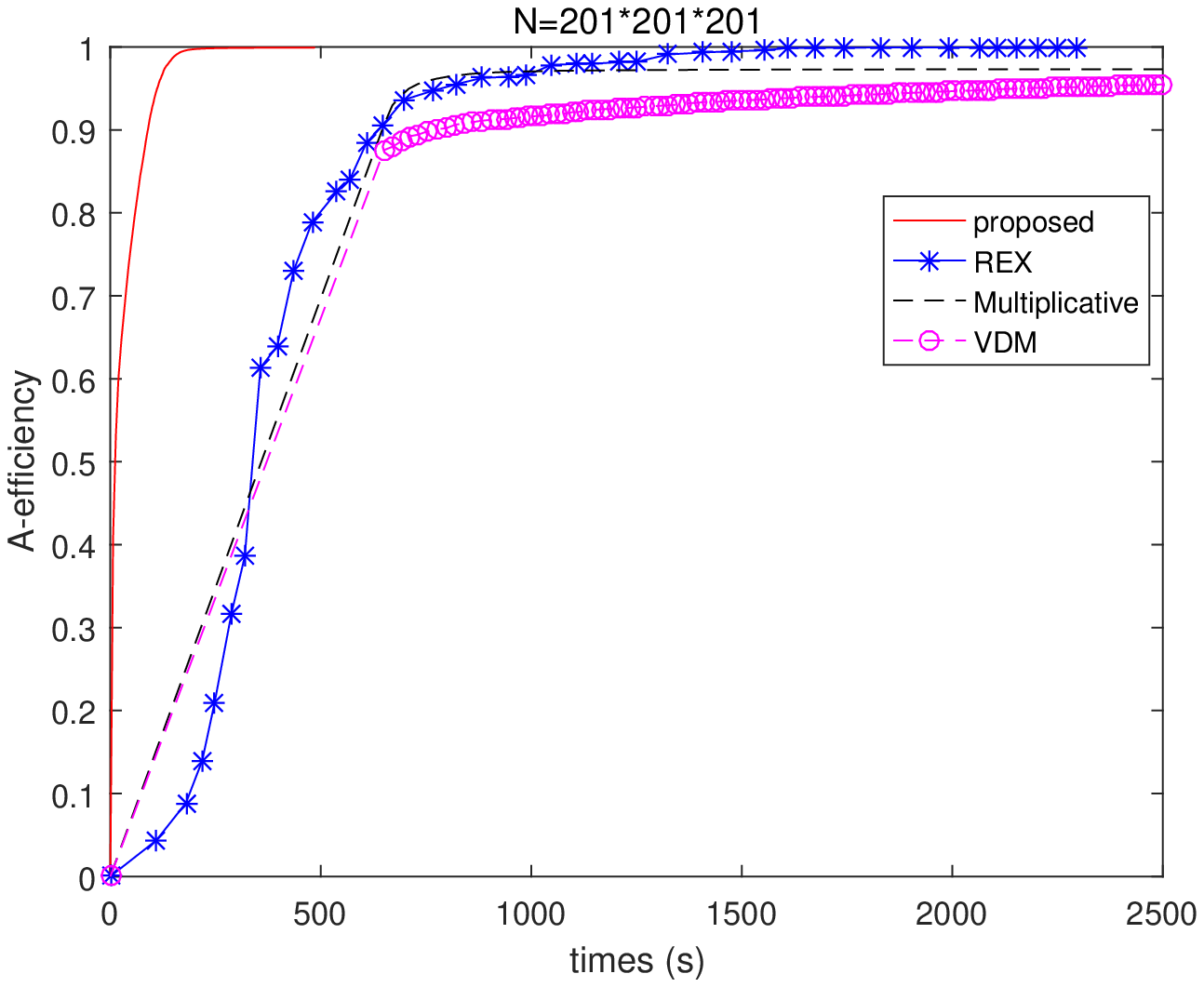}\\
  \caption{The $A$-efficiency for Setting 4 when $N=201\times 201 \times201$. The vertical axis denotes the $A$-efficiency. The horizontal axis corresponds to the computation time (in seconds).}\label{AOP_N201}
\end{figure}

To compare the speed of the proposed algorithm for $D$-optimality and $A$-optimality with the REX, CO, VDM and MUL algorithms, we plot the $D$-efficiency and $A$-efficiency (i.e,
$|\sum\limits_{i=1}^{N}w^{*}_{i}\bff_{i}\bff_{i}^{T}|/|\sum\limits_{i=1}^{N}w_{i}\bff_{i}\bff_{i}^{T}|$ and\\ ${\mbox {tr}} \left( \left[\sum\limits_{i=1}^{N} w^{*}_{i}\bff_{i}\bff^{T}_{i} \right]^{-1}\right)/{\mbox {tr}} \left( \left[\sum\limits_{i=1}^{N} w_{i}\bff_{i}\bff^{T}_{i} \right]^{-1}\right)$, where $w^{*}_i$ $i=1,\cdots,N$ is the theoretical optimal design) versus the time (seconds) for Setting 4 with varying sizes from Figures \ref{DOP_N21} -- \ref{AOP_N201}. From Figures \ref{DOP_N21} -- \ref{DOP_N201}, we can see that all of these algorithms will ultimately converge to the theoretical optimal design (given enough time), but the proposed method is superior to the other methods for large size of $N$ of the design space because of the parallel strategy for the computation as mentioned in Remark 1. However, for a small $N$, the REX algorithm tends to perform better than the proposed method. A similar observation can also be drawn for the $A$-efficiency from Figures \ref{AOP_N21} -- \ref{AOP_N201}.

\begin{figure}[H]
  \centering
  \includegraphics[width=0.6\textwidth]{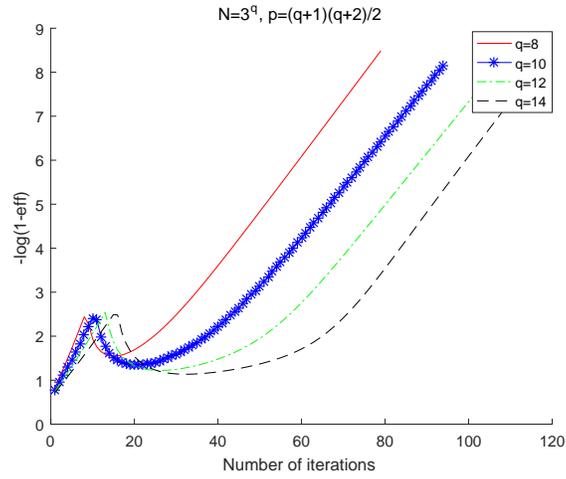}\\
  \caption{The log-efficiency $-\log_{10}(1-eff)$ of the proposed method with $q=8,10,12,14$ for model (\ref{quadratic_regression}), where $eff$ is the lower bound of $D$-efficiency of design. The values $1, 2, 3, \cdots$ on the vertical axis correspond to $D$-efficiency $0.9, 0.99, 0.999, \cdots$. The horizontal axis corresponds to the number of iterations.}\label{log_eff}
\end{figure}

\begin{figure}[H]
  \centering
  \includegraphics[width=0.6\textwidth]{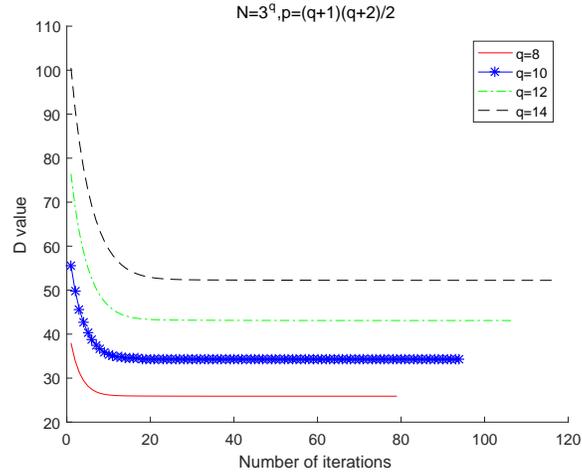}\\
  \caption{The $D$-criterion value of the proposed method with $q=8,10,12,14$ for model in Eq. (\ref{quadratic_regression}). The vertical axis value denotes the $D$-criterion values of the proposed method and the horizontal axis denotes the number of iterations.}\label{D_value}
\end{figure}

To illustrate the performance of the proposed method in the case that the number of factors is large, we consider the full quadratic regression model
\begin{eqnarray}
y(x_1,x_2,\cdots,x_q)=\beta_0+\sum\limits_{i=1}^q \beta_ix_i+\sum\limits_{j=1}^q\sum\limits_{k=j}^q \beta_{j,k}x_jx_k+\varepsilon,
\label{quadratic_regression}
\end{eqnarray}
where $\beta_{1,1},\beta_{1,2},\cdots,\beta_{q,q}$ correspond to the parameters $\beta_{q+1},\beta_{q+2},\cdots,\beta_{p}$, $p=\frac{(q+1)(q+2)}{2}$.
In Figure \ref{log_eff}, the vertical axis is the value of $-\log_{10} (1-eff)$ for $q=8,10,12,14$ and the horizontal axis denotes the number of iterations, where $eff$ represents the lower bound of the $D$-efficiency \citep{Pukelsheim2006}: $D$-efficiency $\geq \frac{p}{\max_{\bx\in \mathcal{X}} \bff^{'}(\bx) \bM(\xi) \bff(\bx)}$, where $\xi$ is the current design.
In other words, the vertical axis in Figure \ref{log_eff}, the values $1, 2, 3, \cdots$ correspond to $D$-efficiency $0.9, 0.99, 0.999, \cdots$.
In Figure \ref{D_value}, the vertical axis is the $D$-criterion values of designs produced by the proposed method for $q=8,10,12,14$ and the horizontal axis denotes the number of iterations.
From Figures \ref{log_eff} and \ref{D_value}, with the increase of iteration times, $-\log_{10} (1-eff)$ gradually increases and $D$-criterion values converge to the $D$-optimal value when the number of factors is large. Thus, the proposed method is still effective even when the number of factors is large.

Based on the numerical evaluations of the five settings considered in this section, we found that the proposed algorithms for $D$-optimality and $A$-optimality converge in all cases, and the optimal design points, as well as the corresponding weights, are close to the theoretical values. Furthermore, in some cases, the proposed method outperforms some existing algorithms for computing approximate  $D$- and $A$-optimal designs. It is noteworthy that the proposed algorithm is simple and it can be implemented without relying on any advanced mathematical programming solvers. Therefore, the proposed algorithms provide a more convenient and effective way to approximate the $D$- and $A$-optimal solutions on the compact design space.

\section{Concluding Remarks}
In this paper, we discuss the approximate optimal design and proposed efficient iterative computational algorithms to obtain the approximate $D$-optimal and $A$-optimal designs for linear models on compact design spaces.
Due to the simplicity and efficiency of the algorithm, the two proposed algorithms are easy to implement. The proposed algorithms are useful tools for many practical applications of optimal design of experiments.

We also provided proof of the monotonic convergence of the proposed algorithm for $D$-optimality and demonstrate that the proposed algorithms provide solution that converges to the optimal design. Furthermore, we prove that the optimal approximate designs converge to the continuous optimal design under certain conditions. A theoretical justification for the convergence of the proposed algorithm for $A$-optimality is not available, but our numerical results strongly support the validity and reliability of the proposed algorithm. These algorithms are implemented in Matlab and the programs are available from the authors upon request. It is worth mentioning that although we focus on $D$-optimal designs and $A$-optimal designs for linear models in this paper, the ideas of the proposed algorithms can be extended to other optimal criteria and other design space in high-dimensional situations.

\section*{Appendix}
\subsection*{Proof of Theorem 1}
Since $\log(|\textbf{A}|)$ is concave in $\textbf{A}$ with $\textbf{A}$ being a positive definite matrix, we have
\begin{eqnarray*}
&&\log \left\{ \left| (1-\lambda)\sum\limits_{i=1}^{N} w^{*}(\bx_{i})\bff(\bx_{i})\bff^{T}(\bx_{i})+\lambda\sum\limits_{i=1}^{N} w(\bx_{i})\bff(\bx_{i})\bff^{T}(\bx_{i}) \right| \right\}\\
&&\geq (1-\lambda) \log\left\{\left|\sum\limits_{i=1}^{N} w^{*}(\bx_{i})\bff(\bx_{i})\bff^{T}(\bx_{i}) \right| \right\}+\lambda \log\left\{ \left|\sum\limits_{i=1}^{N} w(\bx_{i})\bff(\bx_{i})\bff^{T}(\bx_{i}) \right| \right\}.
\end{eqnarray*}
Then, $w^{*}$ is the optimal solution for the $D$-optimal criterion in (2) if and only if
\begin{eqnarray*}
\frac{\log \left\{ \left| (1-\lambda) I_{\mathcal{X}}(\bw^{*}, \bff) + \lambda  I_{\mathcal{X}}(\bw, \bff) \right| \right\} - \log \left\{ \left| I_{\mathcal{X}}(\bw^{*}, \bff) \right| \right\}}{\lambda} \leq 0,
\end{eqnarray*}
where $I_{\mathcal{X}}(\bw, \bff) = \sum\limits_{i=1}^{N} w(\bx_{i})\bff(\bx_{i})\bff^{T}(\bx_{i})$, for all $w(\bx_{i})$ that satisfy $w(\bx_{i})\geq 0, \; i=1,2,\ldots,N$ and $\sum\limits_{i=1}^{N} w(\bx_{i}) =1$, and $\lambda > 0$.
Thus, for $\lambda\downarrow 0$, we have
\begin{eqnarray*}
&   &\lim\limits_{\lambda\downarrow 0} \frac{\log(|(1-\lambda) I_{\mathcal{X}}(\bw^{*}, \bff) + \lambda I_{\mathcal{X}}(\bw, \bff)|) - \log(|I_{\mathcal{X}}(\bw^{*},\bff)|)}{\lambda}\\
& = & \left. \frac{\partial \log \{ |(1-\lambda) I_{\mathcal{X}}(\bw^{*}, \bff) +\lambda I_{\mathcal{X}}(\bw, \bff)| \} }{\partial \lambda} \right|_{\lambda=0}\\
& = & \mbox{tr} \left\{ \left[\sum\limits_{i=1}^{N} w^{*}(\bx_{i})\bff(\bx_{i})\bff^{T}(\bx_{i}) \right]^{-1} \sum\limits_{i=1}^{N} [w(\bx_{i})-w^{*}(\bx_{i})]\bff(\bx_{i})\bff^{T}(\bx_{i}) \right\}\\
& = & \mbox{tr} \left\{ \left[\sum\limits_{i=1}^{N} w^{*}(\bx_{i})\bff(\bx_{i})\bff^{T}(\bx_{i}) \right]^{-1} \sum\limits_{i=1}^{N} w(\bx_{i})\bff(\bx_{i})\bff^{T}(\bx_{i}) \right\} - p \leq 0,
\end{eqnarray*}
which gives the result in Theorem 1.

\subsection*{Proof of Theorem 2}
From Lemma 1, we have
\begin{eqnarray*}
& & \log \left|\sum\limits_{i=1}^{N} w^{(n)}(\bx_{i})\bff(\bx_{i})\bff^{T}(\bx_{i}) \right|-\log \left|\sum\limits_{i=1}^{N} w^{(n-1)}(\bx_{i})\bff(\bx_{i})\bff^{T}(\bx_{i}) \right|\\
& \geq & \sum\limits_{i=1}^{N} w^{(n-1)}(\bx_{i})\mbox{trace}\left\{\bff(\bx_{i})\bff^{T}(\bx_{i})D^{(n-1)}\right\}\log \frac{w^{(n)}(\bx_{i})}{w^{(n-1)}(\bx_{i})}\\
& = & \sum\limits_{i=1}^{N} w^{(n-1)}(\bx_{i})\bff^{T}(\bx_{i})D^{(n-1)}\bff(\bx_{i})\log \frac{w^{(n)}(\bx_{i})}{w^{(n-1)}(\bx_{i})}\\
& = & p\sum\limits_{i=1}^{N} w^{(n)}(\bx_{i})\log \frac{w^{(n)}(\bx)}{w^{(n-1)}(\bx)}\geq 0.
\end{eqnarray*}
Thus, we can conclude that
$$\log \left|\sum\limits_{i=1}^{N} w^{(n)}(\bx_{i})\bff(\bx_{i})\bff^{T}(\bx_{i}) \right|$$
is increasing in $n = 2, 3, \ldots.$ We can get that
\begin{eqnarray*}
\log \left|\sum\limits_{i=1}^{N} w^{(n)}(\bx_{i})\bff(\bx_{i})\bff^{T}(\bx_{i}) \right| \leq \log \left|\sum\limits_{i=1}^{N} \bff(\bx_{i})\bff^{T}(\bx_{i}) \right|.
\end{eqnarray*}
Therefore, under the bounded assumption, the sequence $\log \left|\sum\limits_{i=1}^{N} w^{(n)}(\bx_{i})\bff(\bx_{i})\bff^{T}(\bx_{i}) \right|$ is uniformly bounded and increasing, and hence it is convergent.

Using Lemma 1 and Lemma 2, we can obtain
\begin{eqnarray*}
0 & = & \lim_{n \to \infty}\log\left|\sum\limits_{i=1}^{N} w^{(n)}(\bx_{i})\bff(\bx_{i})\bff^{T}(\bx_{i}) \right| - \log\left|\sum\limits_{i=1}^{N} w^{(n-1)}(\bx_{i})\bff(\bx_{i})\bff^{T}(\bx_{i}) \right|\\
& \geq & p \sum\limits_{i=1}^{N} w^{(n)}(\bx_{i})\log \frac{w^{(n)}(\bx_{i})}{w^{(n-1)}(\bx_{i})}\\
& \geq & \frac{p}{2}[\sum\limits_{i=1}^{N} |w^{(n)}(\bx_{i})- w^{(n-1)}(\bx_{i})|]^{2}\geq 0.
\end{eqnarray*}
Then, we can conclude that
$$
\sum\limits_{i=1}^{N} |w^{(n)}(\bx_{i})-w^{(n-1)}(\bx_{i})|\longrightarrow 0 {\mbox { as }} n \longrightarrow +\infty.
$$

\subsection*{Proof of Theorem 4}
We can check that ${\mbox {tr}}(\textbf{A})^{-1}$ is concave in $\textbf{A}$, where $\textbf{A}$ is positive definite matrix. We have
\begin{eqnarray*}
{\mbox {tr}}[(1-\lambda)I_\mathcal{X}(\bw^{*}, \bff)+\lambda I_\mathcal{X}(\bw, \bff)]^{-1} \leq (1-\lambda){\mbox {tr}}(I_\mathcal{X}^{-1}(\bw^{*}, \bff))+\lambda {\mbox {tr}}(I_\mathcal{X}^{-1}(\bw, \bff)).
\end{eqnarray*}
Then, $w^{*}$ is the optimal solution for the $A$-optimal criterion in Eq. (\ref{AOP}) if and only if for all $w(\bx_{i})(w(\bx_{i})\geq 0\;\mbox{and}\;\sum\limits_{i=1}^{N} w(\bx_{i}) =1)$,
\begin{eqnarray*}
\frac{{\mbox {tr}}[(1-\lambda)I_\mathcal{X}(\bw^{*}, \bff)+\lambda I_\mathcal{X}(\bw, \bff)]^{-1} - {\mbox {tr}}(I_\mathcal{X}^{-1}(\bw^{*}, \bff))}{\lambda} \geq 0
\end{eqnarray*}
for $\lambda>0$.
Thus, for $\lambda\downarrow 0$ we have
\begin{eqnarray*}
&&\lim\limits_{\lambda\downarrow 0} \frac{{\mbox {tr}}[(1-\lambda)I_\mathcal{X}(\bw^{*}, \bff)+\lambda I_\mathcal{X}(\bw, \bff)]^{-1}-{\mbox {tr}}(I_\mathcal{X}^{-1}(\bw^{*}, \bff))}{\lambda}\\
& = & \frac{\partial \{ tr[(1-\lambda)I_\mathcal{X}(\bw^{*}, \bff)+\lambda I_\mathcal{X}(\bw, \bff)]^{-1}\}}{\partial \lambda}|_{\lambda=0}\\
& = & -{\mbox {tr}}[I_\mathcal{X}^{-1}(\bw^{*}, \bff)(I_\mathcal{X}(\bw, \bff)-I_\mathcal{X}(\bw^{*}, \bff))I_\mathcal{X}^{-1}(\bw^{*}, \bff) ]\\
& = & -{\mbox {tr}}( I_\mathcal{X}^{-1}(\bw^{*}, \bff)I_\mathcal{X}(\bw, \bff)I_\mathcal{X}^{-1}(\bw^{*}, \bff)) + {\mbox {tr}}( I_\mathcal{X}^{-1}(\bw^{*}, \bff)) \geq 0,
\end{eqnarray*}
which implies Theorem 4.

\bibliography{Bibliography-MM-MC}

\begin{thebibliography}{}

\bibitem[\protect\citename{Atkinson {\em et~al.}, }2007]{Atkinson}
Atkinson, A~C, Donev, A~N, \& Tobias, R~D. 2007.
\newblock {\em Optimum experimental designs, With SAS}.
\newblock Oxford University Press.

\bibitem[\protect\citename{Box {\em et~al.}, }1978]{Box1978}
Box, G E~P, Hunter, W~G, \& Hunter, J~S. 1978.
\newblock {\em Statistics for experimenters: an introduction to design, data
  analysis, and model building}.
\newblock John Wiley $\&$ Sons.

\bibitem[\protect\citename{Castro {\em et~al.}, }2019]{Castro2017}
Castro, Y~D, Gamboa, F, Henrion, D, Hess, R, \& Lasserre, J~B. 2019.
\newblock Approximate optimal designs for multivariate polynomial regression.
\newblock {\em The Annals of Statistics}, {\bf 47}(1), 127--155.

\bibitem[\protect\citename{Chen, }2003]{Chen}
Chen, Y~H. 2003.
\newblock D-optimal designs for linear and quadratic polynomial models.
\newblock {\em National Sun Yat-Sen University, Taiwan}.

\bibitem[\protect\citename{Cook \& Nachtsheim, }1982]{Cook1982}
Cook, R~D, \& Nachtsheim, C~J. 1982.
\newblock Model robust, linear-optimal designs.
\newblock {\em Technometrics}, {\bf 24}(1), 49--54.

\bibitem[\protect\citename{Dempster {\em et~al.}, }1977]{Dempster1977}
Dempster, A~P, Laird, N~M, \& Rubin, D~B. 1977.
\newblock Maximum likelihood estimation from incomplete data via the EM
  algorithm.
\newblock {\em Journal of the Royal Statistical Society, Series B}, {\bf
  39}(1), 1--38.

\bibitem[\protect\citename{Dette \& Studden, }1997]{Dette}
Dette, H, \& Studden, W~J. 1997.
\newblock {\em The theory of canonical moments with applications in statistics,
  probability, and analysis}.
\newblock Wiley \& Sons.

\bibitem[\protect\citename{Duan {\em et~al.}, }2019]{Duan}
Duan, J~T, Gao, W, \& Ng, H K~T. 2019.
\newblock Efficient Computational Algorithm for Optimal Continuous Experimental
  Designs.
\newblock {\em Journal of Computational and Applied Mathematics}, {\bf 350},
  98--113.

\bibitem[\protect\citename{Fedorov, }1972]{Fedorov}
Fedorov, V. 1972.
\newblock {\em Theory of optimal experiments}.
\newblock Academic Press.

\bibitem[\protect\citename{Gao {\em et~al.}, }2014]{Gao}
Gao, W, Chan, P~S, Ng, H K~T, \& Lu, X. 2014.
\newblock Efficient computational algorithm for optimal allocation in
  regression models.
\newblock {\em Journal of Computational and Applied Mathematics}, {\bf 261}(1),
  118--126.

\bibitem[\protect\citename{Gilmour \& Trinca, }2012]{Gilmour2012}
Gilmour, S.~G., \& Trinca, L.~A. 2012.
\newblock Optimum design of experiments for statistical inference.
\newblock {\em Journal of the Royal Statistical Society: Series C (Applied
  Statistics)}, {\bf 61}(3), 345--401.

\bibitem[\protect\citename{Goos {\em et~al.}, }2016]{Goos2016}
Goos, P, Jones, B, \& Syafitri, U. 2016.
\newblock I-Optimal Design of Mixture Experiments.
\newblock {\em Journal of the American Statistical Association}, {\bf
  111}(514), 899--911.

\bibitem[\protect\citename{Harman {\em et~al.}, }2020]{Harman2018}
Harman, R, Filov{\'a}, L, \& Richt{\'a}rik, P. 2020.
\newblock A Randomized Exchange Algorithm for Computing Optimal Approximate
  Designs of Experiments.
\newblock {\em Journal of the American Statistical Association, Accepted}, {\bf
  115}(529), 348--361.

\bibitem[\protect\citename{Karlin \& Studden, }1966]{Karlin1966a}
Karlin, S, \& Studden, W~J. 1966.
\newblock {\em Tchebycheff systems: With applications in analysis and
  statistics}.
\newblock Interscience Publishers John Wiley $\&$ Sons.

\bibitem[\protect\citename{Kiefer, }1974]{Kiefer1974}
Kiefer, J. 1974.
\newblock General equivalence theory for optimum designs (approximate theory).
\newblock {\em The Annals of Statistics}, {\bf 2}(5), 849--879.

\bibitem[\protect\citename{Kiefer \& Wolfowitz, }1959]{Kiefer1959}
Kiefer, J, \& Wolfowitz, J. 1959.
\newblock Optimum Designs in Regression Problems.
\newblock {\em The Annals of Mathematical Statistic}, {\bf 30}(2), 271--294.

\bibitem[\protect\citename{Kullback, }1967]{Kullback}
Kullback, S. 1967.
\newblock A lower bound for discrimination information terms of variation.
\newblock {\em IEEE Transactions on Information Theory}, {\bf 13}(1), 126--127.

\bibitem[\protect\citename{Meyer \& Nachtsheim, }1995]{Meyer1995}
Meyer, R~K, \& Nachtsheim, C~J. 1995.
\newblock The Coordinate-Exchange Algorithm for Constructing Exact Optimal
  Experimental Designs.
\newblock {\em Technometrics}, {\bf 37}(1), 60--69.

\bibitem[\protect\citename{Pukelsheim, }2006]{Pukelsheim2006}
Pukelsheim, F. 2006.
\newblock {\em Optimal Design of Experiments (Classics in Applied
  Mathematics)}.
\newblock Philadelphia, PA: Society for Industrial and Applied Mathematics.
  [348,349,357].

\bibitem[\protect\citename{Silvey {\em et~al.}, }1978]{Silvey}
Silvey, S~D, Titterington, D~M, \& Torsney, B. 1978.
\newblock An algorithm for optimal designs on a finite design space.
\newblock {\em Communications in Statistics - Theory and Methods}, {\bf 7}(14),
  1379--1389.

\bibitem[\protect\citename{Welch, }1982]{Welch1982}
Welch, W.~J. 1982.
\newblock Algorithmic complexity: three np-hard problems in computational
  statistics.
\newblock {\em Journal of Statistical Computation and Simulation}, {\bf 15}(1),
  17--25.

\bibitem[\protect\citename{Wynn, }1970]{Wynn}
Wynn, H.P. 1970.
\newblock The sequential generation of $D$-optimum experimental designs.
\newblock {\em Annals of Mathematical Statistics}, {\bf 41}(5), 1655--1664.

\bibitem[\protect\citename{Yu, }2011]{Yu}
Yu, Y. 2011.
\newblock D-optimal designs via a cocktail algorithm.
\newblock {\em Statistics and Computing}, {\bf 21}(4), 475--481.

\end{thebibliography}

\begin{table}[H]
\newcommand{\tabincell}[2]{\begin{tabular}{@{}#1@{}}#2\end{tabular}}
 \renewcommand{\tabcolsep}{0.2pc}
\renewcommand{\arraystretch}{0.6}
  \centering
   \caption{$D$-optimal design points and weights, and the values of the $D$-optimality objective function for Setting 1}
\begin{tabular}{|c| c| c| c| c| c|} \hline
              &  Weights &  Weights   & Weights & Weights & Weights  \\
Design points & (REX) & (CO)  & (VDM) & (multiplicative) & (proposed algorithm)   \\  \hline
$(-1, -1)$ & 0.1458 & 0.1458 & 0.1625 & 0.1430 &  0.1457\\
$(1, 1)$   & 0.1458 & 0.1458 & 0.1196 & 0.1520 &  0.1457\\
$(-1, 1)$  & 0.1458 & 0.1458 & 0.1595 & 0.1406 &  0.1457\\
$(1, -1)$  & 0.1458 & 0.1458 & 0.1595 & 0.1436 &  0.1457\\
$(0, -1)$  & 0.0802 & 0.0802 & 0.0798 & 0.0914 &  0.0803\\
$(0, 1)$   & 0.0802 & 0.0802 & 0.1196 & 0.0822 &  0.0803\\
$(-1, 0)$  & 0.0802 & 0.0802 & 0.0798 & 0.1066 &  0.0803\\
$(1, 0)$   & 0.0802 & 0.0802 & 0.0798 & 0.1032 &  0.0803\\
$(0, 0)$   & 0.0962 & 0.0962 & 0.0399 & 0.0373 &  0.0960\\ \hline \hline
$-\log(|\Sigma(\bx) |)$ & 4.4706 & 4.4706 &4.5760   &  4.5512  &  4.4718      \\  \hline
\end{tabular}
\label{set1_DOP}
\end{table}
\begin{table}[H]
\newcommand{\tabincell}[2]{\begin{tabular}{@{}#1@{}}#2\end{tabular}}
 \renewcommand{\tabcolsep}{0.2pc}
\renewcommand{\arraystretch}{0.6}
  \centering
   \caption{$A$-optimal design points and weights, and the values of the $A$-optimality objective function for Setting 1}
\begin{tabular}{|c| c| c| c| c| c|} \hline
              &  Weights    & Weights & Weights & Weights   \\
Design points & (REX)  & (VDM) & (multiplicative)  &  (proposed algorithm)  \\  \hline
$(-1, -1)$ & 0.0940 & 0.1202 & 0.0962  & 0.0939\\
$(1, 1)$   & 0.0940 & 0.1202 & 0.0895  & 0.0939\\
$(-1, 1)$  & 0.0940 & 0.1202 & 0.0953  & 0.0939\\
$(1, -1)$  & 0.0940 & 0.1202 & 0.0978  & 0.0939\\
$(0, -1)$  & 0.0978 & 0.1183 & 0.1181  & 0.0978\\
$(0, 1)$   & 0.0978 & 0.0802 & 0.0972  & 0.0978\\
$(-1, 0)$  & 0.0978 & 0.0802 & 0.0972  & 0.0978\\
$(1, 0)$   & 0.0978 & 0.0802 & 0.0873  & 0.0978\\
$(0, 0)$   & 0.2332 & 0.1603 & 0.2214  & 0.2332\\ \hline \hline
$tr(\Sigma(\bx)^{-1})$ &  17.8922 & 18.7300  & 17.9453   &  17.8922      \\  \hline
\end{tabular}
\label{set1_AOP}
\end{table}

\begin{table}[H]
\newcommand{\tabincell}[2]{\begin{tabular}{@{}#1@{}}#2\end{tabular}}
\renewcommand{\tabcolsep}{0.1pc}
\renewcommand{\arraystretch}{0.1}
\centering
{\scriptsize
\caption{$D$-optimal design points and weights, and the values of the $D$-optimality objective function for Setting 2 }
\begin{tabular}{| c |c c| c c| c c| c c | c c |}\hline
& Design & Weights & Design & Weights & Design & Weights & Design &Weights & Design & Weights    \\
&Points &   (REX) & Points &  (CO)   & Points & (VDM)   & Points &  (MUL) & Points & (Proposed) \\ \hline
&(-0.05, -1.00) & 0.0682 & (-1.00, ~0.06) & 0.1103 & (~0.08, ~0.35) & 0.0403 & (-0.10, ~0.92) & 0.0332 & (-1.00, ~0.00) & 0.0222\\
&(-0.78, ~0.62) & 0.0614 & (-0.94, -0.55) & 0.0694 & (-0.76, -0.64) & 0.0403 & (-0.33, -0.94) & 0.0514 & (-0.96, -0.28) & 0.0214\\
&(~0.37, -0.93) & 0.0840 & (-0.94, ~0.55) & 0.0064 & (~0.42, -0.91) & 0.0403 & (-0.01, ~0.81) & 0.0235 & (-0.96, ~0.28) & 0.0812\\
&(~0.01, ~1.00) & 0.1370 & (-0.73, ~0.22) & 0.0569 & (~0.23, ~0.97) & 0.0403 & (~0.86, ~0.50) & 0.0414 & (-0.80, -0.60) & 0.0591\\
&(~0.92, -0.38) & 0.1153 & (-0.73, -0.10) & 0.0989 & (-0.10, ~0.34) & 0.0403 & (-0.56, ~0.39) & 0.0254 & (-0.60, -0.80) & 0.0591\\
&(~0.86, ~0.51) & 0.1306 & (-0.28, -0.82) & 0.0794 & (-0.35, ~0.94) & 0.0403 & (~0.62, -0.78) & 0.0879 & (-0.60, ~0.80) & 0.0924\\
&(-0.78, -0.63) & 0.1343 & (-0.00, ~0.00) & 0.1667 & (-0.91, -0.40) & 0.0331 & (-0.70, ~0.71) & 0.0499 & (-0.28, -0.96) & 0.0214\\
&(~0.00, ~0.00) & 0.1667 & (~0.30, ~0.32) & 0.0825 & (-0.01, -0.05) & 0.0403 & (-0.96, ~0.27) & 0.0514 & (-0.28, ~0.96) & 0.0073\\
&(-0.96, ~0.29) & 0.0978 & (~0.31, ~0.00) & 0.0061 & (~0.99, -0.13) & 0.0403 & (~0.78, ~0.57) & 0.0375 & (~0.00, -1.00) & 0.0222\\
&(-0.39, ~0.92) & 0.0047 & (~0.51, ~0.34) & 0.1069 & (~0.53, -0.85) & 0.0403 & (-0.87, -0.50) & 0.0517 & (~0.00, ~0.00) & 0.1667\\
&               &        & (~0.51, -0.32) & 0.0237 & (-0.69, ~0.72) & 0.0403 & (~0.04, -0.01) & 0.0448 & (~0.00, ~1.00) & 0.0075\\
&               &        & (~0.57, -0.38) & 0.0738 & (~1.00, -0.01) & 0.0403 & (~0.13, -0.73) & 0.0256 & (~0.28, -0.96) & 0.0812\\
&               &        & (~0.57, -0.78) & 0.0207 & (-0.25, ~0.73) & 0.0403 & (~0.98, ~0.18) & 0.0420 & (~0.28, ~0.96) & 0.1241\\
&               &        &                &        & (-0.36, -0.93) & 0.0806 & (~0.85, -0.52) & 0.0450 & (~0.51, ~0.86) & 0.0014\\
&               &        &                &        & (~0.63, -0.77) & 0.0403 & (-0.99, ~0.12) & 0.0509 & (~0.80, -0.60) & 0.0924\\
&               &        &                &        & (-0.96, -0.28) & 0.0403 & (-0.04, ~0.69) & 0.0213 & (~0.86, ~0.51) & 0.0014\\
&               &        &                &        & (-0.96, ~0.26) & 0.0806 & (~0.01, -0.03) & 0.0449 & (~0.96, -0.28) & 0.0073\\
&               &        &                &        & (-0.30, ~0.95) & 0.0403 & (~0.08, -1.00) & 0.0458 & (~0.96, ~0.28) & 0.1241\\
&               &        &                &        & (~0.60, ~0.80) & 0.0806 & (~0.99, ~0.13) & 0.0422 & (~1.00, ~0.00) & 0.0075\\
&               &        &                &        & (~0.10, -0.45) & 0.0403 & (~0.30, ~0.95) & 0.0427 &                &       \\
&               &        &                &        & (~0.00, ~0.40) & 0.0403 & (-0.05, ~1.00) & 0.0433 &                &       \\
&               &        &                &        & (~0.91, ~0.20) & 0.0403 & (-0.71, -0.71) & 0.0531 &                &       \\
&               &        &                &        & (~0.01, -0.04) & 0.0449 &                &        &                &       \\ \hline \hline
$-\log(|\Sigma(\bx) |)$ & \multicolumn{2}{|c|}{8.2497} & \multicolumn{2}{|c|}{17.9858}  & \multicolumn{2}{|c|}{8.6608} & \multicolumn{2}{|c|}{8.6669} & \multicolumn{2}{|c|}{8.2492}   \\  \hline
\end{tabular}
\label{set2_DOP}}
\end{table}
\begin{table}[H]
\newcommand{\tabincell}[2]{\begin{tabular}{@{}#1@{}}#2\end{tabular}}
\renewcommand{\tabcolsep}{0.1pc}
\renewcommand{\arraystretch}{0.1}
\centering
{\scriptsize
\caption{$A$-optimal design points and weights, and the values of the $A$-optimality objective function for Setting 2 }
\begin{tabular}{| c | c c| c c| c c | c c |}\hline
& Design & Weights & Design & Weights & Design &Weights & Design & Weights    \\
& Points &   (REX) & Points & (VDM)   & Points &  (MUL) & Points & (Proposed) \\ \hline
& (-0.05, -1.00) & 0.0371 & (-1.00, -0.03) & 0.0395 & (-0.55, ~0.15) & 0.0174 & (-0.96, ~0.28) & 0.0856\\
& (-0.78, ~0.62) & 0.0355 & (~0.69, -0.73) & 0.0395 & (-0.37, -0.72) & 0.0305 & (-0.80, -0.60) & 0.0700\\
& (~0.37, -0.93) & 0.0861 & (~0.75, -0.66) & 0.0395 & (~0.78, ~0.62) & 0.0930 & (-0.60, -0.80) & 0.0700\\
& (~0.01, ~1.00) & 0.1053 & (-0.09, ~0.13) & 0.0395 & (-0.59, ~0.27) & 0.0149 & (-0.60, ~0.80) & 0.0878\\
& (-0.89, ~0.46) & 0.0978 & (~0.04, -0.08) & 0.0395 & (-0.29, ~0.96) & 0.0383 & (-0.01, ~0.00) & 0.0013\\
& (~0.81, ~0.59) & 0.0781 & (-0.82, ~0.11) & 0.0395 & (~0.86, ~0.51) & 0.0474 & (~0.00, -0.01) & 0.0013\\
& (~0.92, -0.38) & 0.0984 & (~0.04, -0.02) & 0.0395 & (-0.36, ~0.93) & 0.0398 & (~0.00, ~0.00) & 0.2866\\
& (~0.86, ~0.51) & 0.0374 & (~0.92, ~0.40) & 0.0791 & (-1.00, ~0.03) & 0.0613 & (~0.00, ~0.01) & 0.0013\\
& (-0.78, -0.63) & 0.1324 & (-0.18, ~0.15) & 0.0395 & (~0.05, ~0.00) & 0.0539 & (~0.01, ~0.00) & 0.0013\\
& (~0.00, ~0.00) & 0.2919 & (~0.03, ~0.79) & 0.0395 & (~0.83, -0.55) & 0.0396 & (~0.28, -0.96) & 0.0856\\
&                &        & (~0.37, ~0.88) & 0.0395 & (~0.69, -0.57) & 0.0228 & (~0.28, ~0.96) & 0.0998\\
&                &        & (~0.96, ~0.28) & 0.0395 & (~0.05, -0.01) & 0.1077 & (~0.51, ~0.86) & 0.0108\\
&                &        & (-0.30, -0.95) & 0.0395 & (-0.99, -0.11) & 0.0609 & (~0.80, -0.60) & 0.0878\\
&                &        & (~0.03, -0.04) & 0.0395 & (-0.03, ~1.00) & 0.0354 & (~0.86, ~0.51) & 0.0108\\
&                &        & (~0.06, -0.07) & 0.0395 & (-0.56, ~0.65) & 0.0242 & (~0.96, ~0.28) & 0.0998\\
&                &        & (~0.67, -0.74) & 0.0395 & (~0.67, -0.74) & 0.0376 &         &    \\
&                &        & (~0.32, ~0.95) & 0.0395 & (~0.08, ~0.08) & 0.0532 &         &    \\
&                &        & (-0.35, -0.94) & 0.0395 & (-0.63, ~0.17) & 0.0144 &         &    \\
&                &        & (-0.71, -0.70) & 0.0395 & (~0.08, -0.01) & 0.0538 &         &    \\
&                &        & (~0.03, -0.02) & 0.0514 & (-0.31, -0.95) & 0.0567 &         &   \\
&                &        & (-0.59, ~0.81) & 0.1186 & (-0.39, -0.92) & 0.0584 &         &   \\
&                &        & (-0.95, ~0.00) & 0.0395 & (~0.79, -0.62) & 0.0388 &         &   \\ \hline  \hline
$tr(\Sigma(\bx)^{-1})$&  \multicolumn{2}{|c|}{35.2212 }  & \multicolumn{2}{|c|}{38.1327 }  &  \multicolumn{2}{|c|}{38.0870}  &  \multicolumn{2}{|c|}{35.2207}     \\  \hline
\end{tabular}
\label{set2_AOP}}
\end{table}

\begin{table}[H]
\newcommand{\tabincell}[2]{\begin{tabular}{@{}#1@{}}#2\end{tabular}}
 \renewcommand{\tabcolsep}{0.1pc}
\renewcommand{\arraystretch}{0.1}
\centering
{\scriptsize
\caption{$D$-optimal design points and weights, and the values of the $D$-optimality objective function for Setting 3}
\begin{tabular}{| c |c c| c c| c c| c c | c c |}\hline
& Design & Weights & Design & Weights & Design & Weights & Design &Weights & Design & Weights    \\
&Points &   (REX) & Points &  (CO)   & Points & (VDM)   & Points &  (MUL) & Points & (Proposed) \\ \hline
& (-0.35, -0.35) & 0.1652    & (-0.35, -0.35) & 0.1652 &  (-0.35, -0.35)   &0.1626     & (-0.35, -0.35)   &0.1619        &(-0.35, -0.35)              &0.1627   \\
& (-0.35, ~0.35) & 0.1652    & (-0.35, ~0.35) & 0.1652 &  (-0.35, ~0.35)   &0.1652     & (-0.35, ~0.35)   &0.1648        &(-0.35, ~0.35)              &0.1652  \\
& (~0.12, ~0.12) & 0.0690    & (~0.12, ~0.12) & 0.0690 &  (-0.11, ~0.25)   &0.0002     & (~0.14, ~0.00)   &0.0602        &(~0.12, ~0.12)              &0.0690     \\
& (~0.18, ~0.53) & 0.1396    & (~0.18, ~0.53) & 0.1396 &  (~0.03, -0.30)   &0.0002     & (~0.21, ~0.54)   &0.1475        &(~0.18, ~0.53)              &0.1396     \\
& (~0.35, -0.35) & 0.1652    & (~0.35, -0.35) & 0.1652 &  (~0.11, ~0.11)   &0.0003     & (~0.35, -0.35)   &0.1666        &(~0.35, -0.35)              &0.1652     \\
& (~0.53, ~0.18) & 0.1396    & (~0.53, ~0.18) & 0.1396 &  (~0.11, ~0.12)   &0.0004     & (~0.54, ~0.21)   &0.1417        &(~0.53, ~0.18)              &0.1396                 \\
& (~0.70, ~0.70) & 0.1587    & (~0.70, ~0.70) & 0.1587 &  (~0.12, ~0.12)   &0.0677     & (~0.70, ~0.70)   &0.1573        &(~0.70, ~0.70)              &0.1587                \\
& && & &  (~0.18, ~0.46)   &0.0002     &     &              &     &                 \\
& && & &(~0.18, ~0.53)   &0.1392     &     &              &     &  \\
& && & &(~0.27, ~0.07)   &0.0002     &     &              &     &  \\
& && & &(~0.35, -0.35)   &0.1652     &     &              &     &           \\
& && & &(~0.48, ~0.24)   &0.0002     &     &              &     &          \\
& && & &(~0.53, ~0.18)   &0.1394     &     &              &     &              \\
& && & &(~0.57, ~0.63)   &0.0002     &     &              &     &  \\
& && & &(~0.70, ~0.70)   &0.1585     &     &              &     &\\   \hline   \hline
$-\log(|\Sigma(\bx) |)$ & \multicolumn{2}{|c|}{17.5100} & \multicolumn{2}{|c|}{17.5100}  & \multicolumn{2}{|c|}{17.5121}  &  \multicolumn{2}{|c|}{17.5509}& \multicolumn{2}{|c|}{17.5100}   \\  \hline
\end{tabular}
\label{set3_DOP}}
\end{table}
\begin{table}[H]
\newcommand{\tabincell}[2]{\begin{tabular}{@{}#1@{}}#2\end{tabular}}
\renewcommand{\tabcolsep}{0.1pc}
\renewcommand{\arraystretch}{0.1}
\centering
{\scriptsize
\caption{$A$-optimal design points and weights, and the values of the $A$-optimality objective function for Setting 3}
\begin{tabular}{| c | c c| c c| c c | c c |}\hline
& Design & Weights & Design & Weights & Design &Weights & Design & Weights    \\
& Points &   (REX) & Points & (VDM)   & Points &  (MUL) & Points & (Proposed) \\ \hline
&(-0.35, -0.35) &0.1047 & (-0.35, -0.35)&0.1045  & (-0.35, -0.35) &0.1042   &  (-0.35, -0.35)     &0.1047 \\
&(-0.35, ~0.35) &0.1642 & (-0.35, ~0.35)&0.1641  & (-0.35, ~0.35) &0.1618   &  (-0.35, ~0.35)     &0.1642 \\
&(~0.07, ~0.07) &0.1926 & (-0.17, -0.05)&0.0002  & (~0.08, ~0.05) &0.1924   &  (~0.07, ~0.07)      &0.1926 \\
&(~0.21, ~0.54) &0.1567 & (-0.13, -0.02)&0.0002  & (~0.21, ~0.54) &0.1572   &  (~0.21, ~0.54)      &0.1567 \\
&(~0.35, -0.35) &0.1642 & (-0.12, ~0.43)&0.0002  & (~0.35, -0.35) &0.1670   &  (~0.35, -0.35)      &0.1642 \\
&(~0.54, ~0.21) &0.1567 & (~0.07, ~0.07)&0.1901  & (~0.54, ~0.21) &0.1564   &  (~0.54, ~0.21)      &0.1567 \\
&(~0.70, ~0.70) &0.0609 & (~0.07, ~0.08)&0.0005  & (~0.70, ~0.70) &0.0609   &  (~0.70, ~0.70)      &0.0609 \\
&              &         &    (~0.08, ~0.08)&0.0008 & & & & \\
&              &         &    (~0.09, ~0.09)&0.0002 & & & & \\
&              &         &    (~0.09, ~0.10)&0.0003 & & & & \\
&            &           &    (~0.10, ~0.12)&0.0002 & & & & \\
&            &           &    (~0.13, ~0.16)&0.0002 & & & & \\
&               &       &     (~0.14, -0.23)&0.0002 & & & & \\
&               &       &     (~0.21, ~0.54)&0.1564 & & & & \\
&                   &       & (~0.25, -0.31)&0.0002 & & & & \\
&               &       &     (~0.35, -0.35)&0.1641 & & & & \\
&               &       &     (~0.41, -0.05)&0.0002 & & & & \\
&                   &       & (~0.52, ~0.57)&0.0002 & & & & \\
&               &       &     (~0.54, ~0.21)&0.1564 & & & & \\
&               &       &     (~0.55, ~0.24)&0.0002 & & & & \\
&                   &       & (~0.70, ~0.70)&0.0607 & & & & \\ \hline \hline
$tr(\Sigma(\bx)^{-1})$ & \multicolumn{2}{|c|}{359.1845}    & \multicolumn{2}{|c|}{359.4192}   & \multicolumn{2}{|c|}{359.6661}   & \multicolumn{2}{|c|}{359.1845}       \\  \hline
\end{tabular}
\label{set3_AOP}}
\end{table}

\begin{figure}[H]
  \centering
  \includegraphics[width=0.9\textwidth]{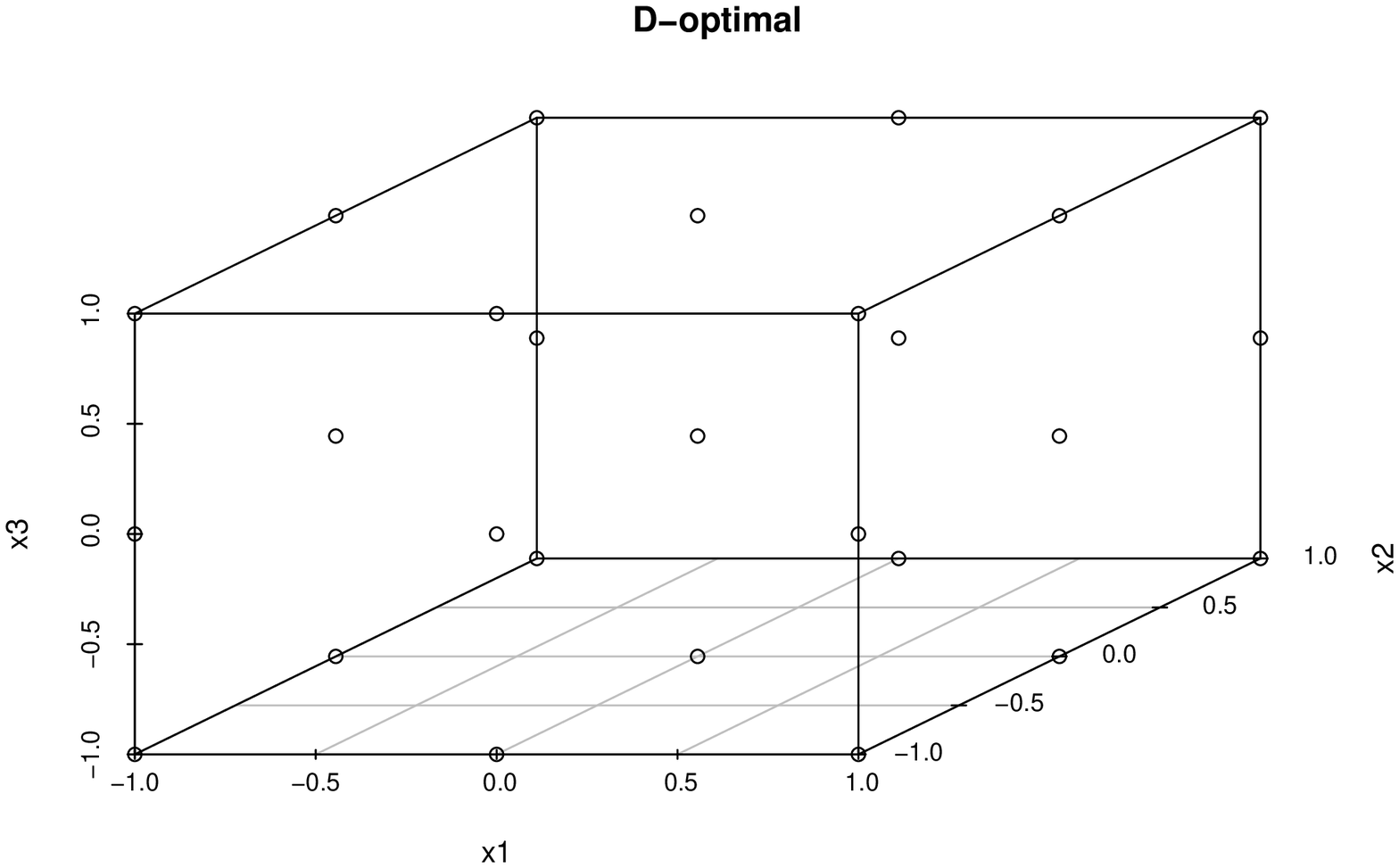}\\
  \caption{$D$-optimal design for Setting 4}\label{set4}
\end{figure}
\begin{table}[H]
\newcommand{\tabincell}[2]{\begin{tabular}{@{}#1@{}}#2\end{tabular}}
 \renewcommand{\tabcolsep}{0.2pc}
\renewcommand{\arraystretch}{0.6}
  \centering
\caption{$D$-optimal design points and weights, and the values of the $D$-optimality objective function for Setting 4}
\begin{tabular}{|c| c| c| c| c| c| }  \hline
   &  Weights  & Weights  & Weights  & Weights   & Weights   \\
Design points & (REX)  & (CO) & (VDM) & (MUL) & (proposed algorithm)  \\  \hline
$(-1, -1, -1)$ &0.0668 & 0.0677 & 0.0489 & 0.0539 & 0.0684\\
$(-1, -1, 1)$ &0.0810  & 0.0613 & 0.0489 & 0.0595 &0.0684\\
$(-1, 1, -1)$ &0.0625  & 0.0698 & 0.0489 & 0.0576 &0.0684\\
$(-1, 1, 1)$ &0.0766   & 0.0635 & 0.0733 & 0.0753 &0.0684\\
$(1, -1, -1)$ &0.0673  & 0.0803 & 0.0489 & 0.0721 &0.0684\\
$(1, -1, 1)$ &0.0815   & 0.0740 & 0.0733 & 0.0726 &0.0684\\
$(1, 1, -1)$ &0.0630   & 0.0825 & 0.0733 & 0.0594 &0.0684\\
$(1, 1, 1)$ &0.0771    & 0.0762 & 0.0733 & 0.0790 &0.0684\\
$(-1, -1, 0)$ &0.0151  & 0.0339 & 0.0489 & 0.0416 &0.0262\\
$(-1, 0, -1)$ &0.0336  & 0.0254 & 0.0244 & 0.0530 &0.0262\\
$(-1, 0, 1)$ &0.0053   & 0.0380 & 0.0489 & 0.0284 &0.0262\\
$(-1, 1, 0)$ &0.0238   & 0.0295 & 0.0244 & 0.0289 &0.0262\\
$(0, -1, -1)$ &0.0288  & 0.0150 & 0.0468 & 0.0533 &0.0262\\
$(0, -1, 1)$ &0.0005   & 0.0276 & 0.0489 & 0.0275 &0.0262\\
$(0, 1, -1)$ &0.0374   & 0.0106 & 0.0244 & 0.0482 &0.0262\\
$(0, 1, 1)$ &0.0091    & 0.0232 & 0.0244 & 0.0137 &0.0262\\
$(1, -1, 0)$ &0.0141   & 0.0086 & 0.0244 & 0.0266 &0.0262\\
$(1, 0, -1)$ &0.0326   & 0.0001 & 0.0244 & 0.0295 &0.0262\\
$(1, 1, 0)$ &0.0228    & 0.0043 & 0.0000 & 0.0270 &0.0262\\
$(1, 0, 1)$ &0.0043    & 0.0127 & 0.0000 & 0.0196 &0.0262\\
$(-1, 0, 0)$ &0.0318   & 0.0073 & 0.0000 & 0.0000 &0.0183\\
$(0, -1, 0)$ &0.0414   & 0.0282 & 0.0000 & 0.0000 &0.0183\\
$(0, 0, -1)$ &0.0045   & 0.0451 & 0.0244 & 0.0000 &0.0183\\
$(0, 0, 1)$ &0.0611    & 0.0199 & 0.0244 & 0.0000 &0.0183\\
$(0, 1, 0)$ &0.0241    & 0.0369 & 0.0489 & 0.0216 &0.0183\\
$(1, 0, 0)$ &0.0339    & 0.0578 & 0.0489 & 0.0164 &0.0183\\
$(0, 0, 0)$ &0.0000    & 0.0005 & 0.0244 & 0.0351 &0.0290\\ \hline \hline
$-\log(|\Sigma(\bx) |)$ & 7.4554 &7.6846 &  8.5046  & 8.0713 &7.4514      \\  \hline
\end{tabular}
\label{set4_DOP}
\end{table}

\begin{table}[H]
\newcommand{\tabincell}[2]{\begin{tabular}{@{}#1@{}}#2\end{tabular}}
 \renewcommand{\tabcolsep}{0.2pc}
\renewcommand{\arraystretch}{0.6}
  \centering
\caption{$A$-optimal design points and weights, and the values of the $A$-optimality objective function for Setting 4}
\begin{tabular}{|c| c| c| c| c| c|} \hline
              &  Weights    & Weights    & Weights  & Weights  \\
Design points & (REX)  & (VDM) & (MUL) & (proposed algorithm) \\  \hline
$(-1, -1, -1)$ &0.0277   &0.0477   &0.0321&0.0402   \\
$(-1, -1, 1)$  &0.0264   &0.0477   &0.0316&0.0402   \\
$(-1, 1, -1)$  &0.0364   &0.0477   &0.0458&0.0402    \\
$(-1, 1, 1)$   &0.0351   &0.0477   &0.0343&0.0402   \\
$(1, -1, -1)$  &0.0450   &0.0238   &0.0268&0.0402    \\
$(1, -1, 1)$   &0.0438   &0.0238   &0.0301&0.0402   \\
$(1, 1, -1)$   &0.0538   &0.0477   &0.0389&0.0402    \\
$(1, 1, 1)$    &0.0525   &0.0477   &0.0424&0.0402   \\
$(-1, -1, 0)$  &0.0521   &0.0238   &0.0465&0.0259   \\
$(-1, 0, -1)$  &0.0421   &0.0238   &0.0282&0.0259   \\
$(-1, 0, 1)$   &0.0447   &0.0000   &0.0307&0.0259  \\
$(-1, 1, 0)$   &0.0347   &0.0000   &0.0271&0.0259   \\
$(0, -1, -1)$  &0.0335   &0.0715   &0.0485&0.0259   \\
$(0, -1, 1)$   &0.0361   &0.0477   &0.0482&0.0259   \\
$(0, 1, -1)$   &0.0161   &0.0238   &0.0529&0.0259  \\
$(0, 1, 1)$    &0.0186   &0.0238   &0.0386&0.0259  \\
$(1, -1, 0)$   &0.0175   &0.0477   &0.0518&0.0259  \\
$(1, 0, -1)$   &0.0075   &0.0477   &0.0455&0.0259   \\
$(1, 1, 0)$    &0.0000   &0.0477   &0.0243&0.0259  \\
$(1, 0, 1)$    &0.0100   &0.0238   &0.0277&0.0259   \\
$(-1, 0, 0)$   &0.0080   &0.0715   &0.0319&0.0430   \\
$(0, -1, 0)$   &0.0253   &0.0238   &0.0000&0.0430   \\
$(0, 0, -1)$   &0.0453   &0.0000   &0.0240&0.0430   \\
$(0, 0, 1)$    &0.0405   &0.0715   &0.0460&0.0430   \\
$(0, 1, 0)$    &0.0602   &0.0477   &0.0221&0.0430   \\
$(1, 0, 0)$    &0.0774   &0.0000   &0.0301&0.0430   \\
$(0, 0, 0)$    &0.1101   &0.0703   &0.0936&0.1096   \\ \hline \hline
$tr(\Sigma(\bx)^{-1})$ & 29.9135 &31.2281 & 31.0049   & 29.9255      \\  \hline
\end{tabular}
\label{set4_AOP}
\end{table}

\begin{table}[H]
\newcommand{\tabincell}[2]{\begin{tabular}{@{}#1@{}}#2\end{tabular}}
 \renewcommand{\tabcolsep}{0.2pc}
\renewcommand{\arraystretch}{0.6}
  \centering
\caption{Values of the $D$- and $A$-optimality objective functions for Setting 5}
\begin{tabular}{|c| c| c| c| c| c|} \hline
Methods & REX & CO  & VDM & multiplicative & proposed algorithm  \\  \hline
$-\log(|\Sigma(\bx) |)$ & 53.2852 & 53.5729  & 52.3201  & 52.7256 & 50.5689    \\  \hline
$tr(\Sigma(\bx)^{-1})$ & 72.0000  & $--$ &  70.4558  & 70.4558 & 70.4139    \\  \hline
\end{tabular}
\label{set5_DEFF}
\end{table}
\end{document}